\documentclass[11pt]{amsart}

\usepackage{fullpage}
\usepackage{url}
\usepackage{amssymb}
\usepackage{cite}

\renewcommand{\a}{\alpha}

\newcommand{\g}{\gamma}
\renewcommand{\d}{\delta}
\newcommand{\D}{\Delta}
\newcommand{\e}{\varepsilon}
\newcommand{\f}{\varphi}
\newcommand{\s}{\sigma}
\renewcommand{\S}{\Sigma}
\renewcommand{\k}{\kappa}
\renewcommand{\l}{\lambda}
 
\renewcommand{\t}{\theta}
\renewcommand{\O}{\Omega}
\renewcommand{\o}{\omega}

\renewcommand{\L}{\Lambda}

\newcommand{\cF}{{\mathcal F}}
\newcommand{\cC}{{\mathcal C}}

\newcommand{\cT}{{\mathcal T}}
\newcommand{\cS}{{\mathcal S}}
\newcommand{\cB}{{\mathcal B}}
\newcommand{\cL}{{\mathcal L}}
\newcommand{\cE}{{\mathcal E}}
\newcommand{\cU}{{\mathcal U}}
\newcommand{\cV}{{\mathcal V}}
\newcommand{\cP}{{\mathcal P}}

\newcommand{\cH}{{\mathcal H}}
\newcommand{\cD}{{\mathcal D}}

\newcommand{\cW}{{\mathcal W}}

\newcommand{\bR}{\mathbb R}

\newcommand{\bZ}{\mathbb Z}

\newcommand{\bE}{\mathbb E}

\newcommand{\be}{\begin{equation}}
\newcommand{\ee}{\end{equation}}

\newcommand{\bel}[1]{\begin{equation}\label{#1}}
\newcommand{\beaa}{\begin{eqnarray*}}
\newcommand{\bea}{\begin{eqnarray}}
\newcommand{\beal}[1]{\begin{eqnarray}\label{#1}}
\newcommand{\bean}{\begin{eqnarray}\nonumber}
\newcommand{\beadl}[1]{\begin{deqarr}\label{#1}}
\newcommand{\eeadl}[1]{\arrlabel{#1}\end{deqarr}}
\newcommand{\eeal}[1]{\label{#1}\end{eqnarray}}
\newcommand{\eead}[1]{\end{deqarr}}
\newcommand{\eea}{\end{eqnarray}}
\newcommand{\eeaa}{\end{eqnarray*}}

\renewcommand{\to}{\rightarrow}

\renewcommand{\phi}{\varphi}
\renewcommand{\epsilon}{\varepsilon}
\renewcommand{\hat}{\widehat}

\newcommand{\<}{\langle}
\renewcommand{\>}{\rangle}

\newcommand{\dm}{{\partial M}}
\newcommand{\dn}{\partial N}
\newcommand{\w}{\widetilde}

\theoremstyle{plain}
\newtheorem{theorem}{Theorem}[section]

\newtheorem{remark}[theorem]{Remark}

\newtheorem{lemma}[theorem]{Lemma}

\newtheorem{proposition}[theorem]{Proposition}
\newtheorem{corollary}[theorem]{Corollary}

\theoremstyle{definition}

\def\endproof{\qed \medskip}
\def\blacksquare{\hbox to .60em {\vrule width .60em height .60em}}

\numberwithin{equation}{section}

\date{\today}

\begin{document}

\title{Static Vacuum Einstein Metrics on Bounded Domains}

\author[ ]{Michael T. Anderson}
\address{Dept.~of Mathematics, Stony Brook University, Stony Brook, NY 11794-3651, USA} 
\email{anderson@math.sunysb.edu}

\thanks{Partially supported by NSF grant DMS 1205947}


\begin{abstract}
We study the existence and uniqueness of solutions to the static vacuum Einstein equations 
in bounded domains, satisfying the Bartnik boundary conditions of prescribed metric and 
mean curvature on the boundary.

\end{abstract}

\maketitle

\section{Introduction}\label{section:intro}

   This paper is a continuation of a series on elliptic boundary value problems for Riemannian 
metrics. Broadly speaking, the issue has both analytic/PDE and geometric/physical aspects. On the 
PDE side, the issue is to identify and solve ``natural" elliptic boundary problems where the unknown 
is a Riemannian metric $g$ on an $(n+1)$-dimensional  manifold $M^{n+1}$ with boundary $\dm$. 
The simplest and most natural problems will involve a quasilinear $2^{\rm nd}$ order (weakly) elliptic 
system for $g$ with suitable boundary data. The fundamental analytic problem is to understand the global 
existence and uniqueness of solutions, (or lack thereof). To be meaningful geometrically and physically, 
the equations and boundary conditions should be diffeomorphism invariant (general covariance). 
Given these constraints, the simplest equations in the interior are the Einstein equations
\be \label{1.1}
Ric_{g} = \l g,
\ee
where $Ric_{g}$ is the Ricci curvature of the metric $g$ and $\l$ is a constant. 
Einstein metrics of course have their origin in general relativity (with metrics of Lorentz signature). 
Moreover, they have long been associated with the problem of finding optimal metrics on a given 
manifold. Any theory describing the existence and uniqueness (or structure of 
the moduli space) of Einstein metrics on closed manifolds is currently far out of reach and the 
perspective here is to see if one can make progress in the presence of a boundary. 

   Both analytically and geometrically, it would appear that the most natural boundary data to 
consider in relation to the Einstein equations \eqref{1.1} are either Dirichlet or Neumann data, i.e.
$$g|_{\dm} = \g_{0} \ \ {\rm or} \ \ A = A_{0},$$
where one prescribes either the induced metric on $\dm$ to be $\g_{0}$ (Dirichlet) or prescribes the 
$2^{\rm nd}$ fundamental form $A$ of $\dm$ in $(M, g)$ to be $A_{0}$ (Neumann). Note this data is 
invariant under diffeomorphisms of $M$ fixing $\dm$. However, such Dirichlet or Neumann data are 
never elliptic for the Einstein equations (in any gauge), cf.~\cite{An2}. Geometrically the most natural 
elliptic boundary data appear to be the mixed Dirichlet-Neumann data
\be \label{1.2}
([\g], H),
\ee
where $[\g]$ is the conformal class of the induced metric on $\dm$ and $H$ is the mean curvature 
of $\dm$ in $M$; cf.~again \cite{An2}. 

\medskip 

  The simplest situation is $3$-dimensions ($n+1 = 3$) where one is considering Einstein metrics on 
3-manifolds with boundary. Here the metrics on the interior are standard, in that they are just of 
constant curvature, so that one is essentially analysing the structure of embedded surfaces in constant 
curvature 3-manifolds with prescribed conformal class and mean curvature. This problem has been 
analysed in detail in \cite{An3}, \cite{An4}, with results that are particularly complete when $M$ 
is a 3-ball with $\dm = S^{2}$. 

   The next simplest case is $n+1 = 4$, so Einstein metrics on $4$-manifolds with 3-manifold boundary. 
This is rather difficult in general, so we focus here on the simplest class of Einstein metrics in dimension 
4, namely the static metrics; these are warped products $N = M\times_{u} S^{1}$ of the form 
\be \label{1.3}
g_{N} = u^{2}d\t^{2} + g,
\ee
with $g = g_{M}$ a Riemannian metric on $M$. 
The Einstein equations \eqref{1.1} with $\l = 0$ on $N$ are then equivalent to the equations
\be \label{1.4}
uRic_{g} = D^{2}u, \ \ \D u = 0,
\ee
for a metric $g$ and potential function $u$ on $M$. The unknowns of the static vacuum Einstein equations 
\eqref{1.4} are thus $(M, g, u)$. Many, but not all, of the results of the paper apply to the case $\l \neq 0$, 
cf.~Remark 3.6 for more details. For simplicity, we concentrate on the situation $\l = 0$. 

\medskip

  Throughout we assume that $N$ is a 4-manifold with boundary $\dn$ and $u > 0$ on $\dn$, 
so that $\dn = S^{1}\times \S$ with $\S$ a closed surface; we assume $N$ and $\dn$ are connected 
and oriented. The metric $g_{N}$ is called {\it strictly static} if 
$$u > 0$$
in the interior of $N$. In this case $N = S^{1}\times M$ topologically, with $M$ a 
3-manifold with connected boundary $\dn = \S$. The simplest case is $N = S^{1}\times B^{3}$, 
with $M = B^{3}$, $\dm = S^{2}$. 
At the same time, we will also consider the situation where the ``horizon"
$$\cH = \{u = 0 \subset N \} \neq \emptyset.$$
Loosely speaking, this situation corresponds to domains surrounding a collection of black holes. 
The static vacuum equations \eqref{1.4} imply $\cH$ is a finite collection of totally geodesic surfaces 
in $M$, giving rise to new boundary components of $M$, but not to $N$. We write $\partial_{tot}M = 
\dm \cup \cH$ with $\dn = S^{1}\times \dm$, so that $\dm$ represents the ``outer" boundary 
of $M$. Again for $\dn = S^{1}\times S^{2}$, the simplest example is the filling $N = D^{2}\times S^{2}$. 
More generally, consider $\dn = S^{1} \times \S$ where $\S$ is any surface. This bounds on the one 
hand $S^{1} \times H$, where $H$ is a 3-dimensional handlebody (corresponding to the strict static 
case), or also $B^{2}\times \S$, (corresponding to the black hole case). In general, the horizon $\cH$ 
could have several components - even of different genus. 

\medskip 

   Let $\cE = \cE^{m,\a} = \cE^{m,\a}(M)$ be the moduli space of $C^{m,\a}$ smooth static vacuum Einstein 
metrics on $M$; this consists of pairs $(g, u)$ on $M$ with $u > 0$ on $M \cup \dm$, cf. Section 2 for the 
precise definition. Throughout the paper we assume $m \geq 2$ and $\a \in (0,1)$. It follows from \cite{An2}, 
cf.~also \cite{AK}, that if $\cE \neq \emptyset$, then $\cE$ is a smooth infinite dimensional Banach manifold. 
This implies that one has a good local structure to the space of solutions of the static vacuum Einstein 
equations. 

  While one could consider elliptic boundary data of the form \eqref{1.2} on $N$, here we focus on 
data introduced by Bartnik \cite{Ba1}, \cite{Ba2}, motivated by issues in general relativity related to 
quasi-local mass. Thus, given a static vacuum solution $(M, g, u)$, consider the boundary data 
\be \label{1.5}
(\g, H),
\ee
where $\g$ is the induced metric and $H$ is the mean curvature on $\dm$. It is shown in \cite{An2}, \cite{AK} 
that the boundary data \eqref{1.5} form elliptic boundary data for the static vacuum equations \eqref{1.4}, 
(in a suitable gauge). Moreover, the boundary map
\be \label{1.6}
\Pi_{B}: \cE^{m,\a} \to Met^{m,\a}(\dm)\times C^{m-1,\a}(\dm),
\ee
$$\Pi_{B}(g, u) = (\g, H),$$
is a smooth Fredholm map of Fredholm index 0. Basic questions of global existence and uniqueness 
of static vacuum metrics with given boundary data $(\g, H)$ are then equivalent to the surjectivity 
and injectivity of $\Pi_{B}$. 

\medskip 

   The key to obtaining such global information on the structure of the map $\Pi$ and the space 
$\cE$ is to understand when the map $\Pi$ is proper (at least when mapping onto suitable domains in 
the target space). This amounts to showing that solutions $(M, g, u)$ are controlled by their boundary 
data $(\g, H)$, i.e.~obtaining apriori estimates for solutions in terms of boundary data. 

   Note that boundary data here are independent of the potential function $u$. Consider first 
the strict static case where $u > 0$ on $N$. The static equations \eqref{1.4} are invariant 
under rescalings $u \to \l u$ of $u$ and 
$$\Pi_{B}(g, \l u) = \Pi_{B}(g, u).$$
Hence $D\Pi_{B}$ always has a non-trivial kernel. Since its Fredholm index is zero, it also has 
everywhere a non-trivial cokernel; cf.~Corollary 2.2  for more details. The image variety 
$$\cV = Im \Pi_{B} \subset \cB \equiv Met^{m,\a}(\dm)\times C^{m-1,\a}(\dm),$$ 
is thus of codimension at least one in the target space $\cB$, and hence $\Pi_{B}$ is never surjective. 
Moreover, it follows that in this situation $\Pi_{B}$ is never proper, since one may always scale 
the length $u$ of the $S^{1}$ fiber so that $u \to 0$ or $u \to \infty$ on $M$. One could 
correct for this by choosing a normalization for $u$, for example $\min u = 1$ or 
$\max u = 1$ and restricting $\cE$ and $\Pi_{B}$ to this subspace. The image of 
$\Pi_{B}$ remains the same, but $\Pi_{B}$ then has Fredholm index $-1$. However to obtain 
useful global information, one needs the Fredholm index to be non-negative. Note that 
this rescaling is not possible in the second (horizon) case, since rescaling creates an edge-type 
singularity along the horizon, cf.~Section 2 for further discussion. 

  There is a another distinct reason that $\Pi_{B}$ is not proper in general. Consider for example 
the space $\cE_{flat}$ of flat static vacuum solutions. Here $(M, g_{flat})$ is a flat 3-manifold and so 
there is an isometric immersion of $M$ (or possibly $\w M$) into $\bR^{3}$ with $\dm$ a 
surface immersed in $\bR^{3}$. For such $(M, g_{flat})$, the function $u$ is any function 
satisfying $D^{2}u = 0$, i.e.~(the pullback of) any affine function of the form $u = a + bz$, for 
$a, b \in \bR$, $z$ a linear function on $\bR^{3}$ of norm 1. The boundary data $(\g, H)$ are 
independent of $u$. However, the space of such $u$ is non-compact; in particular one may have 
sequences of affine functions $u_{i}$ on $M$ with $u_{i} > 0$ for all $i$, but $u = \lim u_{i}$ 
having a zero on $\dm$. Thus, the potential function may degenerate with fixed boundary 
data $(\g, H)$. 

\medskip 

   To overcome these difficulties, we will first show in Proposition 2.6 that, in the strict static case, 
the image $\cV$ of $\Pi_{B}$ is generically transverse to the lines $(\g, \l H)$, $\l \in \bR$, for any 
fixed $(\g, H)$. Hence, it is natural to consider the map to the quotient
\be \label{1.6a}
\Pi_{[B]}: \cE^{m,\a} \to Met^{m,\a}(\dm)\times \cC^{m-1,\a}(\dm),
\ee
$$\Pi_{[B]}(g, u) = (\g, [H]),$$
where $\cC^{m-1,\a}(\dm) = C^{m-1,\a}(\dm)/{\bR}$, with the equivalence relation 
$H_{2} \sim H_{1}$ if and only if $H_{2} = \l H_{1}$ for some $\l$. The map $\Pi_{[B]}$ is still a 
Fredholm map, but now of Fredholm index $+1$. To obtain a map of index 0, 
one may thus append to this map a 1-dimensional space of boundary data and it is most 
natural to choose boundary data depending smoothly on $u$. There is no unique choice for this, 
and the choice we make is to obtain useful information. Thus,  consider 
\be \label{1.7}
\Pi: \cE \to Met(\dm) \times \cC_{+} \times \bR,
\ee
$$\Pi(g, u) = (\g, [H], \mu), \ \ \mu = \int_{\dm}(|d\nu|^{4} + \nu^{2})dv_{\g},$$
where $\nu = \log u$ and $d\nu$ is the differential of $\nu$ on $\dm$. (Again many other 
choices are possible here). Note that the first term in the integral in \eqref{1.7} is invariant 
under rescalings of $u$ but the second is not. 

\medskip 

   The first main result of the paper is the following. Let $\cC_{+}$ be the space of positive 
functions $H: \dm \to \bR^{+}$ modulo the equivalence relation above: two functions 
$H_{i}: \dm \to \bR^{+}$ are equivalent if $H_{2} = \l H_{1}$ for some $\l > 0$. 
Let $\cE_{+}$ denote the space of static vacuum metrics for which the boundary metric has 
positive Gauss curvature and $H > 0$, (modulo rescalings), so that 
$$\cE_{+}^{m,\a} = \Pi^{-1}(Met_{+}^{m,\a}(\dm)\times \cC_{+}^{m-1,\a}(\dm) \times \bR),$$ 
where $Met_{+}$ denotes the space of metrics $\g$ on $\dm$ with positive Gauss curvature: 
$K_{\g} > 0$.

\begin{theorem} The map
\be \label{1.8}
\Pi: \cE_{+} \to Met_{+}\times \cC_{+} \times \bR,
\ee
is a smooth proper Fredholm map of Banach manifolds, of Fredholm index 0.

   This holds in both the strict static case where $\cH = \emptyset$ and the black hole case 
where $\cH \neq \emptyset$. 
   
\end{theorem}

  Theorem 1.1 makes it possible to study the global behavior of the map $\Pi$ by using 
robust topological techniques. Foremost among these is the notion of degree of a map. 
In particular, a smooth, proper Fredholm map $F: X \to Y$ of index 0 between Banach manifolds 
has a well-defined $\bZ_{2}$-valued degree - the Smale degree, \cite{Sm}, defined by 
\be \label{1.9}
deg_{\bZ_{2}} \, F = \# F^{-1}(y), \ \ {\rm (mod \, 2}) ,
\ee
where $y \in Y$ is any regular value of $F$; the properness of $F$ guarantees that the 
cardinality in \eqref{1.9} is finite. This applies to the map $\Pi$ in \eqref{1.8}.

\begin{theorem} In the strict static case, with $M = B^{3}$, $N = S^{1}\times B^{3}$, one has
\be \label{1.10}
deg_{\bZ_{2}} \Pi = 0.
\ee
\end{theorem}

   In fact Theorem 1.2 holds for $M$ of arbitrary topology, e.g.~$M$ equal to a handlebody 
$H$, cf.~Remark 4.6. We conjecture that \eqref{1.10} also holds in the ``Schwarzschild black 
hole case", where $N = S^{2}\times D^{2}$ with $\dn = S^{2}\times S^{1}$. Evidence is given 
for this conjecture in Section 4, but we have not been able to establish the necessary rigidity 
results to prove the conjecture. 

\medskip 

   The techniques used to prove Theorems 1.1 and 1.2 also allow one to use Morse theory and 
Lyusternik-Schnirelman theory methods to analyse the number of Einstein metrics with given 
$\Pi$-boundary values. An example of this is the following rather surprising consequence. 
Let $(M, g, u) = (B^{3}(1), g_{Eucl}, a + bz)$ be a standard flat solution of the static Einstein 
equations \eqref{1.4} on the unit ball in $\bR^{3}$. The induced data on the boundary are 
$(\g, H) = (\g_{+1}, 2)$. 

\begin{theorem} In the strict static case, with $M = B^{3}$, $N = S^{1}\times B^{3}$, 
any boundary data $(\g, [H], \mu)$, $\mu > 0$, near the standard flat data $(\g_{+1}, [2], \mu)$ 
are realized by at least $3$ distinct static vacuum solutions, i.e.
$$\# \Pi^{-1}(\g, [H], \mu) \geq 3,$$
for all data $(\g, [H], \mu)$ near $(\g_{+1}, [2], \mu)$, with $\mu > 0$. Further, for generic 
boundary data $(\g, [H], \mu)$, $\mu > 0$, near standard flat data $(\g_{+1}, [2], \mu)$ there 
are at least $4$ distinct static vacuum solutions, i.e.
$$\# \Pi^{-1}(\g, [H], \mu) \geq 4.$$
 
\end{theorem}

  Theorem 1.3 is related to a recent conjecture of Jauregui \cite{J}; this and more details regarding 
the content of Theorem 1.3 are discussed in Section 5. 

\medskip 

   We conclude the Introduction with a brief summary of the contents of the paper. In Section 2, we 
establish background material and preliminary results for the main results to follow; Propositions 2.5 
and 2.6 are of particular note. Section 3 is devoted to the proof of Theorem 1.1, while Section 4 
computes the degree of $\Pi$ in some standard situations. Section 5 concludes with a discussion of 
Morse theory aspects and the proof of Theorem 1.3. 
   
\section{Preliminary Results}

  In this section, we discuss and prove preliminary results needed for the work to follow. As 
discussed in the Introduction, the static vacuum equations \eqref{1.4} are equations for a pair $(g, u)$, 
consisting of a Riemannian metric $g$ and a scalar function (potential) $u > 0$ on a 3-manifold 
$M$ with boundary $\dm$, possibly with inner boundary components $\cH$ on which $u = 0$.

  To determine the structure of $\cH$, the static equations \eqref{1.4} imply that 
\be \label{2.0}
D^{2}u = 0 \ \ {\rm on} \ \ \cH.
\ee
Evaluating \eqref{2.0} on a pair of tangent vectors $T$ to $\cH$ gives, since $u = 0$ on $\cH$, 
$N(u)A = 0$ on $\cH$, where $A$ is the second fundamental form and $N$ is a unit normal. 
Evaluating \eqref{2.0} on $(N, T)$ where $N$ is normal and $T$ is tangent to $\cH$ gives $TN(u) = 0$, 
so that $N(u) = const$ is a constant on each component of $\cH$ (known as the surface gravity). 
Since $u$ is harmonic, the unique continuation property for harmonic functions implies that $
N(u) \neq 0$ and hence $A = 0$ on $\cH$ so that the horizon is totally geodesic in $(M, g)$.

   The constraint equations, i.e.~Gauss and Gauss-Codazzi equations, at $\dm$ take the form
\be \label{2.1}
|A|^{2} - H^{2} + R_{\g}  = R_{g} - 2Ric(N,N) = 2u^{-1}(\D u + HN(u)),
\ee
\be \label{2.2}
\d(A - H\g) = -Ric(N, \cdot) = -u^{-1}(dN(u) - A(du)).
\ee
To see \eqref{2.1}, one uses $R_{g} = 0$ and $-Ric(N,N) = -u^{-1}NN(u) = u^{-1}(\D_{\dm} u + HN(u))$, 
where the last equation uses $\D_{M}u = 0$. The divergence constraint \eqref{2.2} is a similar 
simple computation. 
   
\medskip 

   Given $M$, let $\bE = \bE^{m,\a}(M)$ be the space of all static vacuum solutions 
$(M, g, u)$ which are $C^{m,\a}$ up to $\dm$. The group ${\rm Diff}_{1}^{m+1,\a}(M)$ of $C^{m+1,\a}$ 
diffeomorphisms of $M$, equal to the identity on $\dm$ acts freely on $\bE^{m,\a}$. The quotient 
$\cE^{m,\a} = \bE^{m,\a} / {\rm Diff}_{1}^{m+1,\a}$ is the moduli space of static vacuum 
Einstein metrics. Observe that the boundary data $(\g, H)$ are well defined on quotient. 

  It is proved in \cite{An2} (or \cite{AK}) that the space $\cE^{m,\a}$ is a smooth infinite dimensional 
Banach manifold (if not empty) and the boundary map 
\be \label{2.3a}
\Pi_{B}: \cE^{m,\a} \to Met^{m,\a}(\dm) \times C^{m-1,\a}(\dm) \equiv \cB^{m,\a},
\ee
$$\Pi_{B}(g) = (\g, H),$$
is a $C^{\infty}$ smooth Fredholm map, of Fredholm index 0. 

\medskip  
 
   The boundary data $(\g, H)$ arise from a natural variational problem. Let $S(h) = -\Delta tr h 
+ \delta \delta(h) - \langle Ric, h \rangle$ be the linearization of the scalar curvature, (cf.~\cite{Be} 
for instance), with adjoint $S^{*}$ given by
$$S^{*}u = D^{2}u -\Delta u \, g - uRic.$$
It is well-known that the static vacuum equations are given by 
$$S^{*}u = 0 \ \ {\rm and} \ \  R = 0.$$ We first discuss the strict static case $\cH = \emptyset$. 

\begin{proposition}\label{p2.1}
On the space of strict static metrics with boundary conditions as in \eqref{1.5}, the Lagrangian
\begin{equation} \label{2.3}
\cL(g, u) = -\int_{M}uR dV_{g}: Met(M) \times C(M) 
\rightarrow \bR,
\end{equation}
has gradient $\nabla \cL$ at $(g, u)$ given by
\begin{equation} \label{2.4}
-\nabla \cL = (S^{*}u + {\tfrac{1}{2}}uRg, \ R, \ uA - N(u)\g,\  2u).
\end{equation}
Thus, if $(h, u')$ is a variation of $(g, u)$ inducing the 
variation $(h^{T}, H'_{h})$ of the boundary data, then 
\begin{equation} \label{2.5}
-d\cL(h, u', h^{T}, H'_{h}) = \int_{M}[\< S^{*}u + {\tfrac{1}{2}}uRg, 
h \> + Ru'] + \int_{\dm}[\< uA - N(u) \g, h^{T} \> + 2uH'_{h}].
\end{equation}
   In particular the static vacuum equations are critical points for $\cL$ with data 
$(\g, H)$ fixed on $\dm$. 
\end{proposition}

{\bf Proof:} The reason for the choice of the minus sign in \eqref{2.3} is to obtain the right 
sign in \eqref{2.17} below. To begin, one has
\begin{equation} \label{2.6}
-D\cL(h, u') = \int_{M}(uR' + u'R + uR(dV)') = \int_{M}uS(h) 
+ {\tfrac{1}{2}}uR \< g, h \> + Ru'.
\end{equation}
A simple integration-by-parts argument gives 
\begin{equation} \label{2.7}
\int_{M}uS(h)  = \int_{M}\< S^{*}u, h \> + 
\int_{\dm}-uN(tr h) - (\d h)(N)u - \< h(N), du \> + tr h N(u).
\end{equation}

  The equations \eqref{2.6} and \eqref{2.7} imply immediately the bulk Euler-Lagrange 
equations - the first two terms in \eqref{2.4}. If the bulk term (over $M$) vanishes, 
then since $u'$ is arbitrary $R = 0$, so this gives 
$$S^{*}u = 0,$$
with $R = 0$, which are the static vacuum equations. 

  For the boundary terms, a standard formula gives:
$$2A'_{h} = \nabla_{N}h + 2A\circ h - 2\d^{*}(h(N)^{T}) - \d^{*}(h_{00}N),$$
where $T$ denotes tangential projection onto the boundary and $h_{00} = h(N,N)$, so 
that (since $H' = tr A' - \< A, h \>$),  
$$2H'_{h} = N(tr h) + 2\delta(h(N)^{T}) - h_{00}H - N(h_{00}).$$
Also by a simple calculation 
$$(\d h)(N) = \d(h(N)^{T}) + \< A, h \> - h_{00}H - N(h_{00}),$$
so that, 
$$2H'_{h} - (\delta h)(N) = N(tr h) + \d(h(N)^{T}) - \< A, h \>.$$
This gives
$$\int_{\dm}u(-N(tr h) - (\d h)(N)) = \int_{\dm}-2uH'_{h} + 
\< du, h(N)^{T}\> - u\< A, h \>.$$
It follows that the boundary term in \eqref{2.7} is given by
$$\int_{\dm}-2uH'_{h} - u\< A, h \> - N(u)h_{00} + N(u)tr h = 
\int_{\dm}-2uH'_{h} - u\< A, h^{T} \> + 
N(u)\< h^{T}, \g \>,$$
which proves the result. 

 {\endproof}

  As noted in the Introduction, in the strict static case $Ker D\Pi_{B}$ is 
always non-empty; the rescalings $u \to \l u$ of $u$ leave the boundary data invariant, 
so that $(0, \l u) \in Ker D\Pi$. Since $D\Pi$ is of Fredholm index 0, it follows that $D\Pi_{B}$ 
always has non-trivial cokernel. 

\begin{corollary}\label{c2.2}
In the strict static case, the image $\cV$ of $\Pi_{B}$ has codimension at least one in the 
target space $\cB$, 
\begin{equation} \label{2.8}
codim(Im \Pi_{B}) \geq 1.
\end{equation}
In fact the tangent space of the image of $\Pi_{B}$ at $\Pi_{B}(g, u)$ is orthogonal to the 
span of the vector $Z = (-uA + N(u)\gamma, -2u) \in T{\mathcal B}$,
\begin{equation} \label{2.9}
(\Pi_{B})_{*}(T\cE)_{(g,u)} = T(Im \Pi_{B}) \perp (-uA + N(u)\gamma, -2u) .
\end{equation}
In particular $Im \Pi_{B}$ has empty interior in ${\mathcal B}$. 
\end{corollary}

{\bf Proof:}  On the moduli space $\cE$ (i.e.~on-shell) one has $\cL = 0$ and hence 
$d\cL = 0$ on variations tangent to $\cE$. The left side of \eqref{2.5} 
thus vanishes, as does the bulk term on the right. Hence, the boundary term on the 
right vanishes, which is exactly the statement \eqref{2.9}. 

{\endproof}

  It would be interesting to know if \eqref{2.9} can be explicitly ``integrated", so that the codimension 1 
restriction for $\cV \subset \cB$ can be described explicitly. 

\medskip 

   Corollary 2.2 does not hold in the black hole case $\cH \neq \emptyset$; in this case, the same 
computations as in the proof of Proposition 2.1 give 
$$\int_{\dm}\langle uA - N(u) \gamma, h^{T} \rangle + 2uH'_{h} = -\int_{\cH}N(u)\<\g, h^{T}\>.$$
Since $N(u) = \k$ is the constant surface gravity, one has
$$\int_{\cH}N(u)\<\g, h^{T}\> = 2N(u)\frac{d}{dt}area(\cH),$$
where the derivative is taken with respect to $g_{t} = g + th$. This term may or may not vanish; it is 
analogous to the variation of a mass term for asymptotically flat metrics. 

  To obtain a variational characterization of the boundary data \eqref{2.3a} in the black hole case, we pass 
to the 4-manifold $N$ with $g_{N} = u^{2}d\theta^{2} + g$ with $\dn = S^{1}\times \dm$. Consider 
first the Einstein-Hilbert action on $Met(N)$ with Gibbons-Hawking-York boundary term
\be \label{2.10}
\cL = -\int_{N}R^{N}dv_{g_{N}} - 2\int_{\dn}H^{N}dv_{\g_{N}},
\ee
where $R^{N}$ and $H^{N}$ are the scalar curvature and mean curvature at $\dn$ of $g_{N}$. 
It is standard that the variation of $\cL$ is given by
\be \label{2.11}
d\cL(h) = \int_{N}\<E, h^{N}\>dv_{g_{N}}  + \int_{\dn}\<\tau, h^{N}\>dv_{\g_{N}},
\ee
where $\tau = A^{N} - H^{N}\g_{N}$ is the conjugate momentum and $E = Ric_{g_{N}} - 
\frac{R_{N}}{2}g_{N}$ is the Einstein tensor of $g_{N}$. One has $H^{N} = N(\nu) + H$, where 
as before, $H$ is the mean curvature of $\dm \subset M$. The matrix corresponding to the form 
$\tau$ is given by
$$\tau = \left(
\begin{array}{cc}
-H & 0 \\
0 & A - (N(\nu) + H)\g 
\end{array}
\right).$$
The derivative \eqref{2.11} does not vanish on the space of deformations $h$ fixing the boundary data $(\g, H)$, 
due to the ``vertical term" (tangent to the $S^{1}$ fiber) at the boundary. To see this, let 
$\o = u d\theta$ be the unit length 1-form tangent to the $S^{1}$ fiber. Then $h^{N} = 2uu'd\theta^{2} 
+ h = 2\frac{u'}{u}\o^{2} + h$. If $h^{T} = 0$ at $\dm$, then $\<\tau, h^{N}\>dv_{\g_{N}} = - 2Hu'dv_{\g} \neq 0$. 

  Thus in place of $\cL$ consider
\be \label{2.12}
\w \cL = -\int_{N}R^{N}dv_{g_{N}} - 2\int_{\dn}H^{N}dv_{\g_{N}} + 2\int_{\dn}Hdv_{\g_{N}}.
\ee
Computing as above then gives
$$d\w \cL(h) = \int_{N}\<E, h^{N}\>dv_{g_{N}} + \int_{\dn}[\<\tau, h^{T}\> + 2H'_{h} + Htr_{\g}h^{T}]dv_{\g_{N}},$$
which does vanish when $h^{T} = H'_{h} = 0$. 

   This gives the analog of Proposition 2.1 in the black hole case. 

\begin{proposition}\label{p2.3}
On the space of black hole static metrics on $N$ with boundary conditions $(\g, H)$, the Lagrangian 
$\w \cL$ has gradient given by 
\begin{equation} \label{2.13}
-\nabla \w \cL = (S^{*}u + {\tfrac{1}{2}}uRg, \ R, \ uA - N(u)\gamma,\  2u).
\end{equation}
Thus, if $(h, u')$ is a variation of $(g, u)$ inducing the 
variation $(h^{T}, H'_{h})$ of the boundary data, then 
\begin{equation} \label{2.14}
-d\w \cL(h, u', h^{T}, H'_{h}) = \int_{M}[\langle S^{*}u + {\tfrac{1}{2}}uRg, 
h \rangle + Ru'] + \int_{\dm}[\langle uA - N(u) \gamma, h^{T} \rangle + 2uH'_{h}].
\end{equation}
   In particular the static vacuum equations are critical points for $\w \cL$ with data 
$(\gamma, H)$ fixed on $\dm$. 
\end{proposition}

{\endproof}

   On the other hand, the analog of Corollary 2.2 does not hold in the black hole case, since 
the functional $\w \cL$ does not vanish on shell, i.e.~is not trivial on the space of solutions $\cE$.

\medskip 

  When $K = Ker D\Pi_{B} = 0$, the map $\Pi_{B}$ is a local diffeomorphism near a given solution $(g, u)$, 
so it has a local inverse map $\cB \to \cE$ defined near $\Pi_{B}(g, u)$. When $K \neq 0$ this is no longer 
the case, but the well-known Lyapunov-Schmidt reduction gives this structure modulo finite dimensional spaces. 
To describe this in the current context, let
\be \label{2.15}
E(g, u) = -(S^{*}u + {\tfrac{1}{2}}uRg, R),
\ee
be the static Einstein operator, so that $\bE = E^{-1}(0,0)$. Of course this may also be written in the 
form
\be \label{2.16}
E(g, u) = (u[Ric_{g} - \frac{R}{2}g - u^{-1}D^{2}u + u^{-1}\D u \, g] , -R).
\ee
In the following it will be more convenient to work on the 4-manifold $(N, g_{N})$, where \eqref{2.16} has 
the simple form
\be \label{2.17}
E(g_{N}) = Ric_{g_{N}} - \frac{R_{N}}{2}g_{N}.
\ee
Let $h^{N}$ be a variation of $g_{N}$ and let 
\be \label{2.19}
L(h^{N}) = E'(h^{N}),
\end{equation}
be the linearization of $E$ at $(N, g_{N}) \in \bE$ with fixed boundary data 
\be \label{2.20}
h^{T} = H'_{h} = 0 \ \ {\rm at} \ \ \dm.
\ee
The operator $L$ gives the $2^{\rm nd}$ variation of $\cL$ in \eqref{2.3} or $\w \cL$ in \eqref{2.12} and 
so is symmetric for this choice of boundary values:
$$\<L(h_{1}), h_{2}\> = \<h_{1}, L(h_{2})\>,$$
for $h_{i}$ satisfying \eqref{2.20}.  Explicitly, $L$ is given by
\be \label{2.21}
2L(h^{N}) = D^{*}Dh^{N} - 2R(h^{N}) - D^{2}tr h^{N} - \d \d h^{N} \, g + \D tr h^{N} \, g - 2\d^{*}\d h^{N},
\ee
where here all geometric quantities are on $(N, g_{N})$. We view $L$ as acting
$$L: S_{0}^{m,\a}(N) \to S^{m-2,\a}(N),$$
where $S^{m,\a}(N) = TMet_{stat}^{m,\a}(N)$ is the tangent space of static metrics on $N$ and 
$0$ denotes the boundary condition \eqref{2.20}. 

    Now write
$$S_{0}^{m,\a} = Im \d^{*} \oplus \cT_{0}^{m,\a},$$
where $\cT^{m,\a}$ denotes the transverse space of divergence-free forms on $N$, $\d_{N}h = 0$, and 
$\cT_{0}^{m,\a}$ denotes the subspace of such forms satisfying \eqref{2.20}. Also $\d^{*}$ acts on $C^{m+1,\a}$ 
vector fields $Z$ on $M$ with $Z = 0$ on $\dm$, (so the boundary condition \eqref{2.20} is preserved). 
This $L^{2}$-orthognal splitting is well-known, cf.~\cite{Be} for instance. 

  The forms $\d^{*}Z$ are in the kernel of $L$, $L(\d^{*}Z) = 0$, and $D^{2}\cL(\d^{*}X, \cdot) = 0$. 
(This corresponds of course to the fact that the functionals $\cL$ and $\w \cL$ are invariant under 
static diffeomorphisms in ${\rm Diff}_{1}(N)$). Moreover, by the symmetry property, 
$$Im L \perp Im \d^{*}.$$
To see this, one has $\<L(h), \d^{*}Z\> = \int_{N}\<\d(L(h)), Z\>$, since 
$Z = 0$ on $\dm$. But $\d L(h) = \d(E'_{h}) = (\d E)'_{h} - \d'_{h}E = 0$. Hence
\be \label{2.22}
L: \cT_{0}^{m,\a} \to \cT^{m-2,\a}.
\ee
This operator is elliptic, so Fredholm, and symmetric. 

   The {\it nullity} of $(M, g, u) = (N, g_{N})$ with respect to $\Pi_{B}$ is the dimension of the 
kernel $K = Ker D\Pi_{B}$, consisting of forms $\k = (k, 2\nu')$, $\nu' = \frac{u'}{u}$, satisfying
$$L(\k) = 0, \ \ \d_{N}\k = 0, \ \ {\rm with} \ \ k^{T} = 0, \ H'_{k} = 0 \ {\rm at} \ \dm.$$
A regular point of $\Pi_{B}$ has nullity zero. A set $\S \subset \cE$ is a {\em non-degenerate critical manifold} 
if $\S$ is a finite dimensional manifold and each point in $\S$ is a critical point of $\Pi_{B}$ with nullity equal 
to $dim \S$. 

  Let $K^{\perp}$ be the $L^{2}$ orthogonal complement of $K \subset \cT^{m-2,\a}$ 
and let $K_{0}^{\perp}$ be the $L^{2}$ orthogonal complement of $K \subset \cT_{0}^{m,\a}$. 
The symmetry (self-adjointness) property above implies that 
\be \label{2.23}
L: K_{0}^{\perp} \to K^{\perp},
\ee
is an isomorphism. 

   Now suppose $K \neq 0$ at a given solution $(g_{0}, u_{0})$ with boundary data $(\g_{0}, H_{0})$. Given 
boundary data $(\g, H)$ near $(\g_{0}, H_{0})$ choose an extension to a metric $g_{\g,H}$ on $N$ for 
instance as follows: $g_{N} = u_{0}^{2}d\theta^{2} + g_{M}$ with 
$$g_{M} = \eta(dt^{2} + \f^{2}\g) + (1 - \eta)g_{0},$$
where $\f = 1$ and $N(\f) = \frac{1}{2}H$ at $t = 0$ so that $g_{M}$ realizes the boundary data $(\g, H)$. 
Here $t$ is the geodesic normal coordinate for $g_{0}$ so that $g_{0} = dt^{2} + \g_{t}$ and $\eta$ is a 
bump function with $\eta = 1$ near $\dm$. This gives a map $(\g, H) \to g_{\g, H} \in 
Met^{m,\a}(M)\times C^{m,\a}(M)$. Now define an operator 
$$P: \cB \times K \times K_{0}^{\perp} \to K^{\perp},$$
$$P(\g, H, \k, \ell) = \pi_{K^{\perp}}[E(g_{\g,H} + \k + \ell)],$$
where $\pi_{K^{\perp}}$ is the orthogonal projection to $K^{\perp}$, $\k = (k, 2\nu') \in K$ and 
$\ell = (h, 2\nu') \in K_{0}^{\perp}$. The derivative $(D_{3}P)_{\g_{0}, H_{0},0,0}$ in the third 
component is an isomorphism by \eqref{2.23}. Hence by the implicit function theorem, there is a 
(smooth) map 
$$I_{0}: \cB \times K \to K_{0}^{\perp},$$
such that $I_{0}(\g,H, \k)$ is the unique solution (locally) to 
$$P(\g, H, \k, I_{0}(\g, H, \k)) = 0,$$
so that 
$$E(g_{\g,H} + \k + I_{0}(\g, H, \k)) \in K.$$
Set 
$$I: \cB \times K \to Met^{m,\a}(M) \times C^{m,\a}(M),$$
$$I(\g, H, \k) = g_{\g,H} + \k + I_{0}(\g, H, \k)$$
and set 
\be \label{2.24}
z: \cB \times K \to K,
\ee
$$z(\g, H, \k) = E(I(\g, H, \k)) \in K.$$ 

  From this we have the following result on the structure of $\cE$ near $(g_{0}, u_{0})$. 
  
\begin{proposition}\label{p.2.4}
For $(\g, H)$ near $(\g_{0}, H_{0})$ as above, one has $I(\g, H, \k) \in \cE$ if and only if $z(\g, H, \k) = 0$. 
The locus $z^{-1}(0)$ is a smooth submanifold of $\cB \times K$ of codimension $k = dim K$ with $0\times K 
\subset T(z^{-1}(0))$. Further, for arbitrary fixed boundary data $(\g, H)$, $z^{-1}(0) \cap \Pi^{-1}(\g, H)$ 
consists exactly of the critical points of $\cL$ (in case $\cH = \emptyset$) or $\w \cL$ (in case 
$\cH \neq \emptyset$). 
\end{proposition}

{\bf Proof:}  The first statement is immediate from the discussion above. The second statement is a 
consequence of the result \eqref{2.3a} that $\cE$ is a smooth manifold; namely the divergence-gauged 
Einstein operator $Met_{stat}^{m,\a}(N) \to S^{m-2,\a}(N)$ is a submersion at any Einstein metric, 
cf.~\cite{An2}, \cite{AK}. The last statement is immediate from Propositions 2.1 and 2.3. 

{\endproof}

  Abusing notation slightly, let $\cE^{nd}$ be the set of ''non-degenerate" static Einstein metrics, i.e.~the set of 
such metrics $(g, u)$ such that $Ker D\Pi_{B} = (0, \l u)$ in the strict static case, and $Ker D\Pi_{B} = (0, 0)$
in the black hole case. In the former case one can of course remove the scaling freedom of $u$ by restricting 
to the subspace of $\cE$ for which, for instance, $\int_{M}u = 1$; however we will not carry this out 
in practice. Thus regular points of $\Pi_{B}$ are defined to be those with nullity equal to 1 in the 
strict static case. Let $\cE^{d} = \cE \setminus \cE^{nd}$ be the complement consisting of degenerate 
metrics. Let $\cV^{d}$ and $\cV^{nd}$ denote the image of these sets under $\Pi_{B}$. The next result shows 
that non-degenerate metrics (i.e.~the regular points of $\Pi_{B}$) are open and dense in $\cE$; more precisely the 
degenerate metrics are contained in a submanifold of codimension 1 in $\cE$. 

\begin{proposition}\label{p2.5}
   In the strict static case $\cH = \emptyset$, $Im(\Pi_{B})$ is of codimension one in the target space $\cB$. 
The "singular points" of $\Pi_{B}$, i.e.~ the locus $\cE^{d}$, is of codimension at least one in $\cE$. At any point 
in $\cV^{nd}$, one has
\be \label{2.25}
T\cV^{nd} = Z^{\perp} \subset T\cB,
\ee
where $Z$ is given by \eqref{2.9}. 

  In the black hole case $\cH \neq \emptyset$, $Im(\Pi_{B})$ is of codimension zero in the target space $\cB$. 
Again, the singular points of $\Pi_{B}$, i.e.~ the locus $\cE^{d}$, is of codimension at least one in $\cE$. 
At any point in $\cV^{nd}$, one has
\be \label{2.26}
T\cV^{nd} = T\cB.
\ee
\end{proposition}

{\bf Proof:} The proof is based on the classical idea that if $(M, g, u)$ is degenerate, so that $L$ in \eqref{2.21} 
has a zero eigenvalue ($K = Ker L \neq 0$) then perturbing $M$ slightly by moving the boundary 
$\dm$ inward by the normal exponential map should give a non-degenerate solution with $K = 0$. 
We follow some of the methods used in \cite{W2} to carry this out in this context.

   The issue is local, so fix any degenerate solution $(M, g_{0}, u_{0})$ with $dim K = k > 0$. Then as above, 
solutions in $\cE$ near $(M, g_{0}, u_{0})$ are parametrized by $(\g, H, \k)$ and 
$$I(\g, H, \k) \in \cE$$
exactly when $z(\g, H, \k) = 0$. Consider the 2-parameter variation
$$I(\g_{s}, H_{s}, \k_{s}+t\k) = g_{s,t},$$
where $I(\g_{s}, H_{s}, \k_{s})$ is the curve in $\cE$ given by moving in along unit inward 
normal $-N$ a geodesic distance $s$, so this curve is in $z^{-1}(0)$. The curve $I(\g_{0}, H_{0}, t\k)$
is tangent to $\cE$ but may not be in $\cE$ at second order. 

  Let $\cS_{k} = \{(\g, H, \k) \in z^{-1}(0): dim K = k\}$. Observe that 
$$\cS_{k} = \{(\g, H, \k): z(\g, H, \k) = 0, \ \ D_{2}z(\g, H, \k) = 0\},$$
where $D_{2}$ is the derivative in the $\k$ direction, i.e.~the full derivative of $z$ in directions of $K$ 
vanishes. Note that $D_{2}z = L$, for $L$ as in \eqref{2.21}, acting on forms satisfying the boundary conditions 
\eqref{2.20}. 

   At this point, we divide the discussion into two cases. 

\medskip 
   
  Case I: $D_{2}^{2}z(\g_{0}, H_{0}, 0) \neq 0$. 

In this case, there is some direction $\k_{0}$ so 
that the $2^{\rm nd}$ derivative of $z$ in the $\k_{0}$ direction is non-zero. Thus the first derivative 
of $\partial_{\k_{0}}z(\g_{0}, H_{0}, \k)$ at $\k = 0$ is non-zero. But $(0, \k_{0})$ is tangent to $z^{-1}(0)$
and hence the set of $(\g, H, \k)$ with $I(\g, H, \k)$ near $(g_{0}, u_{0})$ such that $z(\g, H, \k) = 
\partial_{k_{0}}z(\g, H, \k) = 0$ lies in a submanifold of $z^{-1}(0)$. Since this set contains $\cS_{k}$ 
near $(g_{0}, u_{0})$, it follows that $\cS_{k}$ has codimension at least one. 

\medskip 

   Case II: $D_{2}^{2}z(\g_{0}, H_{0}, 0) = 0$.

  In this situation, we construct a functional on $\cE$ which vanishes at the critical locus $\cE^{d}$  
of $\Pi_{B}$ but with non-zero gradient. This again gives the codimension one property.  
We work first in strict static case and follow that with the black hole case. 

   To begin, note that the action $\cL$ in \eqref{2.3} vanishes identically on $\cE$, so gives no information. 
Instead, consider 
\be \label{2.27}
\cF = \cL - \int_{\dm}Hdv_{\g} = -\int_{M}uRdv_{g} - \int_{\dm}Hdv_{\g}.
\ee
This is non-zero on $\cE$ and the variation onshell (i.e.~on $\cE$) is given by 
$$\d \cF(h) = -(\int_{\dm}Hdv_{\g})' = -\int_{\dm}H'_{h} + {\tfrac{1}{2}}Htr_{\g}hdv_{\g}.$$ 
Observe this vanishes on $K = Ker D\Pi_{B}$, so that critical points of $\Pi_{B}$ are critical points of 
$\cF$. More importantly, at $g \in \cE$, a form $\k \in K$ exactly when the second variation 
$D^{2}\cF(\k, \cdot) = 0$, when $\cF$ is viewed as a functional on the space of metrics with fixed 
boundary metric and mean curvature, as in \eqref{2.20}. Thus $K$ gives the kernel of the second 
variation of $\cF$. 

  Now compute the variation of $\cF$ on the family $g_{s,t}$ in the direction $s$ at $s = 0$, for 
$t$ small. The variation of the metric $g$ is given by $h = -A$ while the variation of the 
potential $u$ is given by $-N(u)$. This gives  
\be \label{2.29}
\frac{d}{ds}\cF(g_{s,t})|_{s=0} = \int_{M}\<E(g, u), (h,u')\> - \int_{\dm}\<uA - N(u)\g, h^{T}\> + 2uH'_{h}
\ee
$$ - \int_{\dm}H'_{h} + {\tfrac{1}{2}}Htr_{\g}hdv_{\g}.$$
Consider the second order behavior of \eqref{2.29} in the direction $t$. The bulk term 
over $M$ vanishes at first order since $E'_{\k} = L(\k) = 0$, as does the first boundary term. Similarly, 
these terms also vanish at second order, exactly by the Case II condition, which implies that 
$D^{2}E(\k, \k) = 0$. Hence 
\be \label{2.30}
\frac{d}{ds}\cF(g_{s,t})|_{s=0} = -{\tfrac{1}{2}}\int_{\dm}N(H) + H^{2} + o(t^{2}),
\ee
since $H'_{h} = -\frac{1}{2}N(H)$. We now compute the second derivative of \eqref{2.30} in the direction 
$\k \in K$, i.e.~in the $t$ direction. First, by definition, one has $2A = \cL_{N}g$, so that 
$A'_{k} = \cL_{N}k$ (in the normal geodesic gauge) so that $A''_{k} = \cL_{N}k' = 0$. 
Similarly $H''_{k} = 0$. 

  Next, by the standard normal Riccati equation, $N(H) = -|A|^{2} - Ric(N,N)$ and by the constraint (Gauss) 
equation $|A|^{2} - H^{2} + R_{\g} = R_{g} - 2Ric(N,N)$. So $-Ric(N,N) = \frac{1}{2}(|A|^{2} - H^{2} 
+ R_{\g} - R_{g})$, giving $N(H) = \frac{1}{2}(-|A|^{2} - H^{2} + R_{\g} - R_{g})$. Note that 
$R_{\g}'' = 0$ since $k^{T} = 0$ at $\dm$ while $R_{g}'' = 0$, since $\k$ is tangent to $\cE$ 
to second order. It follows that 
$$\frac{d^{2}}{dt^{2}}\frac{d}{ds}\cF(g_{s,t})|_{s=t=0} = \frac{1}{4}\int_{\dm}|A'_{k}|^{2}dv_{\g} \geq 0.$$
Suppose first that 
\be \label{2.31}
\frac{d}{ds}\frac{d^{2}}{dt^{2}}\cF(g_{s,t})|_{s=t=0} = \frac{1}{4}\int_{\dm}|A'_{k}|^{2}dv_{\g} > 0,
\ee
for some $\k \in K$. Then the function $\frac{d^{2}}{dt^{2}}\cF(g_{s,t})|_{t=0}$ has non-zero derivative 
in the $s$ or inward direction so that the set of metrics in $\cE$ near the background $(g_{0}, u_{0})$ where 
$\frac{d^{2}}{dt^{2}}\cF(g_{s,t})|_{t=0} = D^{2}\cF(\k,\k) = 0$ is of codimension 1 in $\cE$. However, to first 
order in $s$, $K$ also represents the kernel of $D^{2}\cF$ at $g_{s}$. Namely, if $\k_{s} \in K_{s}$, then 
$L_{s}(\k_{s}) = 0$ so that $L'_{\k}(\k) + L(\k') = 0$. But $L'_{\k}(\k) = D^{2}E(\k, \k) = 0$ by the Case II 
property, so that $L(\k') = 0$, which proves the claim. It then follows that $\cS_{k}$ also has (at least) 
codimension 1. 

   Suppose then on the other hand that 
\be \label{2.32}
k^{T} = (A'_{k})^{T} = 0,
\ee
for all $\k \in K$, so that the Cauchy data for $\k$ vanish at $\dm$. The main issue is to show that 
\eqref{2.32} leads to $\k = 0$, (and hence a contradiction). To do this, we need to obtain the same information 
regarding the Cauchy data for $u'$. Given any $(M, g, u) \in \cE$, let $\w K$ be the space of 
linearized Einstein deformations satisfying \eqref{2.32} at $\dm$, so that in general $\w K \subset K$ but in the 
context in which we are working, $\w K = K$ at the given solution $(M, g_{0}, u_{0})$. Consider then 
the modified functional 
$$\w \cF = -\int_{M}uRdv_{g} - \int_{\dm}Hdv_{\g} - \int_{\dm}\<A, \nabla_{N}A\>dv_{\g}.$$
We view $\w \cF$ as a functional on the space of pairs $(g, u)$ with fixed metric $\g$ and second 
fundamental form $A$ at $\dm$, analogous to but stronger than \eqref{2.20}. As before, the kernel of the 
second variation $D^{2}\w \cF$ on this space is exactly $\w K$ at $(M, g_{0}, u_{0})$. 

  Now note that $-N\<A, \nabla_{N}A\> = -|\nabla_{N}A|^{2} - \<A, \nabla_{N}\nabla_{N}A\>$. As before 
$A''_{k} = 0$ and $\nabla_{N}\nabla_{N}A''_{k} = 0$. Carrying out the same process as above and using \eqref{2.32} 
gives then 
$$\frac{d^{2}}{dt^{2}}\frac{d}{ds}\w \cF(g_{s,t})|_{s=t=0} = \frac{1}{2}\int_{\dm}|A'_{k}|^{2}dv_{\g}
+ \int_{\dm}|\nabla_{N}A'_{k}|^{2}dv_{\g} = \int_{\dm}|\nabla_{N}A'_{k}|^{2}dv_{\g}.$$
If this term is positive, then the same arguments as above prove the result. Suppose then instead
$$\nabla_{N}A'_{k} = 0$$
holds, together with \eqref{2.32}. Taking the derivative of the normal Riccati equation (cf.~\cite{Pe} for 
instance) $\nabla_{N}A + A^{2} + R_{N} = 0$ in the direction $k$ gives $R'_{N} = 0$, and hence by the 
static vacuum equations (and constraint equations on $\dm$) one obtains
\be \label{2.33}
Ric'_{k} = 0
\ee
at $\dm$. 

   Now on $M$, one has  
$$uRic' + u'Ric = (D^{2})'u + D^{2}u'.$$
Since $(D^{2})' = 0$ at $\dm$ by \eqref{2.32}, this together with \eqref{2.33} gives
\be \label{2.34}
u'Ric = D^{2}u'.
\ee
at $\dm$. Thus $u'$ is another static solution, with the same underlying metric $g$, at $\dm$. 
It follows from a result of Tod, \cite{To}, that generically (off a codimension 1 subset of static metrics) the 
potential $u$ is uniquely determined, up to constants, by the metric $g$. More precisely, let $d$ 
denote the covariant exterior derivative acting on $T(M)$-valued 1-forms. Then $dRic$ is the 
well-known Cotton tensor (a vector valued 2-form on $M$). Taking $d$ of the static vacuum 
equations gives 
$$(dRic)_{ij} = (Ric_{ii} - Ric_{jj})d\nu_{k}$$
where $e_{i}$ is a local orthonormal basis diagonalizing the Ricci curvature $Ric$, cf.~[30, Eqn.~(12)]. 
Hence if the eigenvalues $\rho_{i}$ of $Ric$ are distinct, the potential $\nu$ is uniquely determined 
up to a constant by $Ric_{g}$ and $dRic_{g}$. The static metrics for which at least two eigenvalues 
of $Ric$ are equal are of infinite codimension in the space of all static metrics (cf.~[30, Eqn.~(34)]. 

   Thus, without loss of generality, we may assume from \eqref{2.34} and the discussion above that 
$d\log u' = d\log u$ at $\dm$. This gives $u' = cu$ and $N(u') = cN(u)$ at $\dm$. Now to any 
variation $\k = (k, u')$ one may subtract a variation $\eta_{c} = (0, cu)$ rescaling the length of the 
fiber $S^{1}$ so that $\hat \k = \k - \eta_{c} = (k, u' - cu)$. It follows then that for any variation 
$\hat \k$, one has 
\be \label{2.35}
k^{T} = (A'_{k})^{T} = u' = N(u') = 0,
\ee
at $\dm$. Thus, all of the Cauchy data of the variation $\hat \k$ vanish at $\dm$. It follows then 
from the unique continuation theorem in \cite{AH} that 
$$\hat \k = 0,$$
on $M$. In turn, this implies that the original variation $\k$ is trivial, i.e.~of the form $\k = (0, cu)$. 
This completes the proof of Proposition 2.5 in the strict static case. 

   The proof in the black hole case $\cH \neq \emptyset$ is basically the same. One replaces 
$\cL$ in \eqref{2.27} by $\w \cL$. All of the arguments above then carry through, and prove in this 
case also that either $\cS_{k}$ is of codimension (at least) one or \eqref{2.35} holds. It follows again 
that $\k = (0, cu)$ for some $c$. However, in this case, because of the presence of the horizon 
$\cH$, one must have $c = 0$ so that $\k = 0$.

{\endproof}

   Let $\cD^{m,\a} = Met^{m,\a}(\dm)\times \cC^{m-1,\a}(\dm)$, where as in \eqref{1.6a} the 
second term is the space of equivalence classes $[H] = [\l H]$. Let $\pi: \cB \to \cD$ 
denote the projection, so that the fibers of $\pi$ are given by $\pi^{-1}(\g, [H]) = \{(\g, \l H)\}$ 
for $\l \in \bR^{+}$. Thus $\Pi_{[B]} = \pi \circ \Pi_{B}$. Let $\w \cE = \cE / \sim$, where the 
equivalence relation is given by rescalings of $u$, so $u \sim \l u$ for $\l > 0$. This action of 
$\bR^{+}$ on $\cE$ is free, so that $\w \cE$ is also a smooth Banach manifold.

\begin{proposition}\label{p2.6}
In the strict static case ($\cH = \emptyset$) the image space $\cV = Im \Pi_{B}$ is transverse to the 
fibers of $\pi$, except possibly at the critical value locus of $\Pi_{B}$. Thus the map
\be \label{2.36}
\Pi_{[B]}: \w \cE^{nd} \to \cD,
\ee
is a local diffeomorphism.
\end{proposition}

{\bf Proof:}  By Proposition 2.5, it suffices to show that 
\be \label{2.37}
T\cV \oplus T[H] = T\cB,
\ee
away from the critical values in $\cV$. This is equivalent to showing that
$T[H]$ is not contained in $T\cV$. Tangent vectors to $T[H]$ are of the form 
$(0, \l H)$ and by \eqref{2.9} these are never tangent to $\cV = (\Pi_{B})_{*}(T\cE^{nd} )$, since 
$$\int_{\dm}\<uA - N(u)\g, h^{T}\> + 2uH'_{h} = \int_{\dm}2u\l H > 0,$$
since $u > 0$ and $H > 0$. This proves the result. 

 {\endproof}

   Proposition 2.6 is of course false in the black hole case, since $D\Pi_{B}$ maps generically onto 
$T\cB$. A consequence of Proposition 2.6 is the following relation between the regular points of 
$\Pi_{B}$ and $\Pi$ (in the strict static case). 

\begin{corollary} \label{c2.7} 
In the strict static case, a point $(M, g, u) \in \cE$ is a regular point of $\Pi$ if and only 
if $(M, g, u)$ is a regular point of $\Pi_{B}$, (i.e.~a point of nullity $1$). More generally, for 
any $(M, g, u) \in \cE$, 
\be \label{2.38}
{\rm nullity} D\Pi_{B} = {\rm nullity} D\Pi + 1.
\ee
\end{corollary}

{\bf Proof:} If $\k \in Ker D\Pi$, then $k^{T} = 0$, $H'_{k} = \l H$ for some $\l$, 
and $\nu'_{\k} = 0$. The transversality property \eqref{2.37} implies that $H'_{k} = 0$ and hence 
$\k \in Ker D\Pi_{B}$. Thus 
$$K_{\Pi} \subset K_{\Pi_{B}}.$$
Conversely, if $\k \in K_{\Pi_{B}}$, then $\k \in K_{\Pi}$ if and only if $\mu'_{\k} = 0$. This can always be 
achieved by rescaling $u$ ($u \to \l u$), so changing $\k = (k, u') \to (\k, u' - \l)$, for some $\l$, cf.~also 
Remark 4.6.  

{\endproof}

\section{Properness of $\Pi$} 

  To obtain a deeper understanding of the possible boundary values $(\g, H)$ of static 
vacuum solutions, one needs to know that the map $\Pi$ (or $\Pi_{B}$) is proper, or determine 
domains in the target space onto which it is proper. Only with this property can one obtain 
effective information on the global behavior of $\Pi$.

   Proving that $\Pi_{B}$ is proper amounts to showing that a static vacuum metric $(M, g, u)$ can 
be controlled in terms of its boundary data $(\g, H)$, (which amounts to control of the ``inverse 
map" to $\Pi_{B}$). As discussed below, this is relatively easy in the interior of $M$, away from $\dm$ 
(especially in the case $\cH = \emptyset$) but it is more difficult to control the structure of 
$(M, g, u)$ near $\dm$ in terms of $(\g, H)$. The main issue is to obtain some control or 
information about the behavior of the potential function $u$ near $\dm$ in terms of $(\g, H)$. 

  As discussed in the Introduction, $\Pi_{B}$ is in fact never proper, at least in the strict static case, 
due to the possible degeneration of $u$ with fixed or controlled boundary data $(\g, H)$. As 
motivated there, we instead consider the closely related map
\be \label{3.1}
\Pi: \cE \to \cD \times \bR = Met(\dm) \times \cC_{+} \times \bR,
\ee
$$\Pi(g, u) = (\g, [H], \mu),$$
where $\mu = \int_{\dm}(|d\nu|^{4} + \nu^{2})$, $\nu = \log u$. 
The map $\Pi$ is a smooth Fredholm map, of Fredholm index 0. 

   In the following, we prove Theorem 1.1, that $\Pi$ is proper at least on the domain $\cE^{+}$ of 
static vacuum solutions for which $K_{\g} > 0$ and $H > 0$, i.e.~the boundary has positive 
Gauss and mean curvature. 

\medskip 

{\bf Proof of Theorem 1.1.}

Suppose $(M, g_{i}, u_{i})$ is a sequence of static vacuum solutions with boundary data 
\be \label{3.2}
(\g_{i}, [H_{i}], \mu_{i}) \to (\g, [H], \mu),
\ee
in the $C^{m,\a}\times C^{m-1,\a}\times \bR$ topology on $\dm$. One then needs to 
prove that $(M, g_{i}, u_{i}) \to (M, g, u)$ in a subsequence, in the $C^{m,\a}$ topology on $M$. 

   The proofs in the strict static and black hole cases are somewhat different. We thus 
first prove the result in the strict static case $\cH = \emptyset$; following that, 
we discuss the case black hole case $\cH \neq \emptyset$.

\medskip 

Case A. $\cH = \emptyset$.

      Since the metrics $\g_{i}$ are uniformly controlled, by the Sobolev embedding theorem on $\dm$, 
a bound on 
\be \label{3.3}
\mu_{i} = \int_{\dm}(|d\nu_{i}|^{4} + \nu_{i}^{2})dv_{\g_{i}}
\ee
implies a uniform bound on $\sup |\nu_{i}|$ and hence uniform pointwise control on $\{u_{i}\}$ 
away from $0$ and $\infty$ on $\dm$. By the maximum principle for harmonic functions, the same bound 
extends to all of $M$. 

  The main point is then to obtain a uniform bound on the full curvature $Rm$ of $\{g_{i}\}$. To begin, 
write $H_{i} = \l_{i}H_{i}'$, where $|H_{i}'|_{C^{m-1,\a}(\dm)} = 1$. Then either $\l_{i} \to 0$, 
$\l_{i} \to \l_{0} > 0$ or $\l_{i} \to \infty$ (in a subsequence). We analyze these situations case-by-case. 
   
   (I). $\l_{i} \to \l_{0} > 0$. 
   
      In this case both the metrics $\g_{i}$ and mean curvatures $H_{i}$ are uniformly controlled 
along the sequence. 

\begin{proposition}\label{p3.1}
In the strict static case, there  is a constant $\L = \L(\g, H, \mu)$ such that
\be\label{3.4}
|Rm| \leq \L,
\ee
on $M$. In addition, at the boundary $\dm$, 
\be\label{3.5}
|A| \leq \L, \ \ |N(u)| \leq \L, \ \ and \ \ d_{foc} \geq \L^{-1},
\ee
where $d_{foc}$ is the distance to the focal locus of the inward normal exponential map 
from $\dm$ into $M$. The estimates \eqref{3.4}, \eqref{3.5} also hold for higher derivatives 
of $Rm$, $A$ and $N(u)$ up to order $m-2$, $m-1$ respectively.
\end{proposition}

{\bf Proof:} For points $x \in M$ a bounded distance away from $\dm$, the estimate 
\eqref{3.4} follows immediately from the apriori interior estimates in \cite{An1} 
which states that  
\begin{equation}\label{3.6}
|Rm|(x) \leq \frac{k}{t^{2}(x)}, \ \ |d\log u|(x) \leq \frac{k}{t(x)},
\end{equation}
where $t(x) = dist(x, \dm)$ and $k$ is an absolute constant. The estimates \eqref{3.6} 
are scale invariant, and hold also for all higher derivatives of $Rm$ and $\log u$. Consider 
then the behavior at $\dm$. By the scalar constraint equation \eqref{2.1}, a bound on 
$|Rm|$ at $\dm$ implies a bound on$|A|$ at $\dm$, given control of $(\g, H)$. 
Similarly, by standard comparison geometry, a bound on $|Rm|$ on $M$ gives a 
lower bound on the distance $d_{foc}$ to the focal locus of the normal exponential 
map $exp_{\partial M}$. 

   It remains then to prove the curvature bound \eqref{3.4} at or arbitrarily near the boundary 
$\dm$. (The higher derivative estimates then follow by standard elliptic regularity theory).
We prove \eqref{3.4} by a blow-up argument. If the curvature bound in \eqref{3.4} is
false, then there is a sequence $(M, g_{i}, u_{i}, x_{i}) \in \cE$ with
controlled boundary data $(\g_{i}, H_{i})$ such that
$$|Rm_{g_{i}}|\to \infty.$$
Choose points $x_{i} \in M$ such that $|Rm_{g_{i}}|(x_{i} = \max |Rm_{g_{i}}| = \l_{i}^{2} \to \infty$. 
Consider the rescaled metrics $g_{i}' = \lambda_{i}^{2}g_{i}$; then 
\be\label{3.7}
|Rm_{g_{i}'}|(x_{i}) = 1, \ \ {\rm and} \ \ |Rm_{g_{i}'}|(y_{i}) \leq 1,
\ee
for any $y_{i} \in (M, g_{i}')$. By standard scaling properties, $H_{g_{i}'} = \l_{i}^{-1}H_{g_{i}}$, 
$R_{\g_{i}'} = \l_{i}^{-2}R_{\g_{i}}$ and $t_{i}' = \l_{i}t_{i}$. Observe that $t_{i}'(x_{i}) \leq \sqrt{k}$ 
by \eqref{3.6}, so that $x_{i}$ remains within a uniformly bounded 
distance to the boundary $\dm$ with respect to $g_{i}'$.  

  The sequence $(M, g_{i}', u_{i})$ is now uniformly controlled, in that its curvature is uniformly 
bounded, the potential functions $u_{i}$ are uniformly bounded away from $0$ and $\infty$. In 
addition, the boundary geometry consisting of the boundary metric, $2^{\rm nd}$ fundamental 
form and conjugacy radius of the normal exponential map are also uniformly controlled. 
It is possible that within a ball $B_{x_{i}}(R)$ of fixed radius $R$ there are a number of 
components of the boundary, and moreover, some (or all) of these boundary components 
could converge to a single component as $i \to \infty$ (of higher multiplicity). For the 
moment, we assume this is not the case, so that the distance to the cut locus of the 
normal exponential map from $\dm$ is uniformly bounded below. 

    It follows from the convergence theorem in \cite{AT} for manifolds-with-boundary that, for any 
$\a' < \a < 1$, a subsequence converges in $C^{1,\a'}$ to a $C^{1,\a}$ static limit $(X, g, u, x)$ 
with boundary $(\partial X, \g, u)$, cf.~also \cite{K}. The convergence is uniform on compact subsets. 
By the normalization in \eqref{3.7}, the limit $(X, g)$ is complete and without singularities up to the 
boundary $\partial X$. One has $\partial X = \bR^{2}$, the boundary metric $\g$ is flat, $H = 0$,
so $\partial X$ is a minimal surface in $X$. The limit harmonic potential $u$ satisfies $0 < c_{0} < u < 
c_{0}^{-1}$, for some $c_{0}$. The bound \eqref{3.7} and the static equations imply that $u$ 
extends at least $C^{1,\a}$ up to $\partial X$. These remarks, and those below, apply to 
each component of $\partial X$ if there is more than one. 

  Moreover, elliptic regularity associated with the Einstein equation implies that the convergence to the 
limit is in $C^{m,\a}$, since the boundary data $(\g, H)$ are elliptic boundary data and the convergence 
of $(\g, H)$ to the limit is in $C^{m,\a}\times C^{m-1,\a}$. In particular, it follows from \eqref{3.7} 
that 
\be\label{3.8}
|Rm|(x) = 1,
\ee
where $x = \lim x_{i}$ and $Rm = Rm_{X}$.

  On the blow-up limit $(X, g)$, the Gauss equation \eqref{2.1} holds and becomes
$$0 \leq {\tfrac{1}{2}}u|A|^{2} = \D_{\partial X}u,$$
on $\partial X$. Hence $u$ is a $C^{m,\a}$ smooth bounded subharmonic function on 
$\bR^{2}$; it is well-known \cite{GT} that the only such functions are constant, $u = const$. 
It follows also that $A = 0$, so that $\partial X$ is totally geodesic. 

  Next, from the divergence constraint \eqref{2.2}, we also now have $0 = \d (A - H\g) = 
 -u^{-1}D^{2}u(N, \cdot)$, so that $0 = D^{2}u(N, \cdot) = dN(u) - A(du) = dN(u)$. Thus 
$N(u) = const$. In particular, it follows from the static vacuum equations that the 
ambient curvature $Rm$ satisfies $Rm = 0$ at $\partial X$. 

  One sees then that the full Cauchy data $(\g, u, A, N(u))$ for the static vacuum equations
is fixed and trivial: $\g$ is the flat metric, $A = 0$ and $u$, $N(u)$ are constant.
Observe that this data is realized by the family of flat metrics on $(\bR^{3})^{+}$
with either $u = const$ or $u$ equal to an affine function on (a domain in) $(\bR^{3})^{+}$.
By ~cite{AH}, two Einstein metrics on a manifold with boundary with equal Cauchy data (locally) 
are locally isometric. It follows that the limit $(X, g, u)$ is flat near the boundary. 
However, this contradicts \eqref{3.8}. 

   To conclude the proof, we need to remove the assumption of a lower bound on the 
distance cut locus $d_{i}$ to the normal exponential map. Since as noted above $|R| \to 0$ 
at each component of $\partial X$ and \eqref{3.7} holds, if $d_{i} \to 0$ it follows that 
$|\nabla R| \to \infty$ on the sequence $(M, g_{i}', u_{i})$. In this case, we rescale the metrics 
further so that $d_{i} = 1$. It follows that $R \to 0$ with $|\nabla R|(y_{i})| = 1$ at some base 
points $y_{i}$. This is a contradiction to elliptic regularity of the static vacuum equations and so 
cannot occur. 

{\endproof}
 
  (II). $\l_{i} \to 0$. 
  
    In this case, one has $H_{i} \to 0$. If $\sup |R_{i}| \to \infty$, then exactly the same arguments as in 
Case (I) give a contradiction, so one must have $|R_{i}|$ is bounded. From the scalar constraint \eqref{2.1} 
at $\dm$, one has 
\be \label{3.9}
{\tfrac{1}{2}}u(|A|^{2} - H^{2} + R_{\g}) = \D u + HN(u),
\ee
on $\dm$. Now since the curvature $|Rm_{g_{i}}|$ of $g_{i}$ is uniformly bounded, $N(u_{i})$ is also uniformly 
bounded, while $H_{i}$ becomes arbitrarily small so that $HN(u) \sim 0$ for $i$ large. Also $|A|^{2} - H^{2} 
+ R_{\g} \sim |A|^{2} + R_{\g}$. Since $R_{\g_{i}}$ is uniformly bounded away from zero, one obtains 
$\D u > 0$ on $\dm$ for $i$ sufficiently large. This contradicts the maximum principle and hence 
Case (II) cannot occur, i.e.~$H_{i}$ is uniformly bounded away from 0.  

   (III). $\l_{i} \to \infty$. 

    Since $R_{\g} > 0$, a basic result of Shi-Tam [ST] on the positivity of the Brown-York mass applies, 
so that 
\be \label{3.10}
\int_{\dm}H dv_{\g} \leq \int_{\dm}H_{0} dv_{\g},
\ee
where $H_{0}$ is the mean curvature of the unique isometric embedding of $(\dm, \g)$ into $\bR^{3}$. So 
for a given $\g$, this gives an upper bound on $\l$, so that Case (III) also cannot occur. 

\medskip 

   The analysis above shows that only Case (I) can occur. Proposition 3.1 gives a uniform bound on the 
curvature $Rm_{g_{i}}$ of the metrics $g_{i}$ up to $\dm$, and uniform control of $u_{i}$. To obtain 
compactness, one still needs to prove that the manifold-with-boundary structure does not degenerate, 
i.e.~that the normal exponential map from $\dm$ has injectivity radius bounded below,
$$inj_{\dm} \geq i_{0}.$$
More precisely, let $-N$ be the inward unit normal to $\dm$ in $M$ and consider the associated normal 
exponential map to $\dm$, $-tN \to \exp_{p}(-tN)$, giving the geodesic normal to $\dm$ at 
$p$. This is defined for $t$ small, and let $\tau(p)$ be the maximal time interval on which 
$\exp_{p}(-tN) \in M$ is length-minimizing, (so that in particular the geodesic does not hit $\dm$ 
again before time $\tau(p)$). Thus, $\tau: \dm \to \bR^{+}$. 

\begin{lemma}\label{l3.2}
There is a constant $t_{0}$, depending only on $\L$ in \eqref{3.4} and a positive lower bound $H_{0}$ 
for $H$, such that 
\be \label{3.11}
\tau(p) \geq t_{0}.
\ee
\end{lemma}

{\bf Proof:} First, the bound \eqref{3.4} gives a lower bound, say $d_{0}$, on the focal radius $d_{foc}$ 
of the normal exponential map. Suppose then 
$$\min \tau < d_{0}.$$
If the minimum is achieved at $p$, then the normal geodesic $\s$ to $\dm$ at $p$ intersects 
$\dm$ orthogonally again at a point $p'$. Let $\ell = \tau(p)$ be the length of $\s$. The 
$2^{\rm nd}$ variational formula of energy gives
\be \label{3.12}
E''(V,V) = \int_{0}^{\ell}(|\nabla_{T}V|^{2} - \langle R(T,V)V,T\rangle) dt - 
\langle \nabla_{V}T, V \rangle|_{0}^{\ell} \geq 0,
\ee
where $T = \dot \s$ and $V$ is any variation vector field along $\s$ orthogonal to $\s$. 
Choose $V = V_{i}$ to be parallel vector fields $e_{i}$, running over an orthonormal basis at 
$T_{p}(\dm)$. The first term in \eqref{3.12} then vanishes, while the second sums to 
$-Ric(T,T) \leq \L$. The boundary terms sum to $\pm H$, at $p$ and $p'$. Taking into 
account that $T$ points into $M$ at $p$ while it points out of $M$ at $p'$, this gives
$$0 \leq \L \ell - (H(p) + H(p')).$$
Since $H$ is bounded away from $0$, $H \geq H_{0}$, this gives a contradiction if $\ell$ 
is sufficiently small. This proves the estimate \eqref{3.11}. 

{\endproof}

  To complete the proof, one needs an upper bound on the diameter of $(M, g_{i})$ and a lower 
volume bound. The lower volume bound follows trivially from control obtained near $\dm$ 
from Proposition 3.1 and Lemma 3.2. 

\begin{lemma}\label{l3.3}
There is a constant $D$, depending only on a lower bound for $H$, such that 
\be \label{3.13}
diam_{g}(M) \leq D.
\ee
\end{lemma}

{\bf Proof:} The proof of \eqref{3.13} is by contradiction; it would be of interest to find a 
direct proof. If \eqref{3.13} is false, there is a sequence of static vacuum solutions $(M, g_{i}, u_{i})$ 
with boundary data $(\g_{i}, H_{i}, \mu_{i})$ converging to a limit $(\g, H, \mu)$ in the target space topology. 
It follows from the results above that a subsequence converges in $C^{m,\a}$ to a complete 
non-compact static vacuum limit $(\hat M, g, u)$ (of infinite diameter) with boundary data 
$(\g, H)$, $H > 0$. Since $\hat M$ is non-compact, it has a nontrivial end $E$ (possibly many). 
The potential function $u$ on $\hat M$ satisfies $0 < c_{0} \leq u \leq C_{0}$ for some fixed 
constants $c_{0}$, $C_{0}$. 

   We first claim that there is a sequence $r_{j} \to \infty$ such that 
\be \label{3.14}
area(S(r_{j})) \to \infty \ \ {\rm as} \ \ j \to \infty,
\ee
where $S(r) = \{x \in \hat M: dist_{g}(x, \dm) = r\}$. This basically follows from the structure 
of static vacuum ends described in \cite{An1.5}, to which we refer for some further details. 
To start, since $u \geq c_{0}$, by \cite{An1.5} an end $E$ of $\hat M$ is either asymptotically 
flat in the usual sense, or ``small" in the sense that 
$$\int_{1}^{\infty}(area S(r))^{-1}dr = \infty.$$
Clearly \eqref{3.14} holds if the end is asymptotically flat. If the end is small, consider the 
annuli $A(r, kr) = \{x\in \hat M: s(x) \in (r, kr)\}$ ($k$ fixed) rescaled to size 1, so with 
respect to the rescaled metric $g_{r} = r^{-2}g$. It is proved in \cite{An1.5} that when lifted 
to suitable covers (unwrapping of short circles or tori), the metrics $g_{r}$ converge in a 
subsequence to a Weyl solution of the static vacuum equations on $A(1, k)$, i.e. a solution 
on the manifold $I\times S^{1}\times S^{1}$ with metric invariant under rotations about the 
second or both $S^{1}$ factors. The first case corresponds to rank 1 collapse (in the sense of 
Cheeger-Gromov) along circles, while the second case corresponds to rank 2 collapse along 
tori. The second case is only possible for the so-called Kasner metrics (cf.~[5, Ex.2.11]) which 
have unbounded potential $u$.  Hence toral collapse is ruled out and the collapse is along circles. 
In this case, the length of the first $S^{1}$ factor is bounded below with respect to $g_{r}$ so 
grows linearly with respect to $g$, while the length of the second $S^{1}$ factor is bounded 
below with respect to $g$. This gives \eqref{3.14}.

   The surfaces $S(r_{j})$ for any fixed $j$ embed in $(M, g_{i})$ for $i$ sufficiently large and 
are of course homologous to $\dm$. Since $H > 0$ at $\dm$, $\dm$ serves as a barrier for 
the existence of stable (in fact minimizing) surfaces in $(M, g_{i})$. Thus, one may choose $j$ 
sufficiently large, and then $i$ also sufficiently large and find a surface $\S \subset A(0, r_{j}) 
\subset (M, g_{i})$ which minimizes area among surfaces homologous to $\dm$. However, by 
a result of Galloway \cite{Ga} (using the $2^{\rm nd}$ variational formula for area) there 
are no such area minimizing surfaces $\S \subset (M, g_{i})$. This gives a contradiction and shows 
that \eqref{3.13} must hold.

{\endproof}

   Given the results above, it follows from the convergence theorem in \cite{AT}, cf.~also \cite{K}, that 
a subsequence of $(M, g_{i}, u_{i})$ converges, modulo diffeomorphisms, to a limit $(M, g, u) \in \cE$ 
in the $C^{1,\a}$ topology. The convergence in $Met^{m,\a}(M) \times C^{m,\a}(M)$ then follows 
from elliptic regularity, cf.\cite{GT}, \cite{Mo} for instance. 

\medskip 

  Case (B). $\cH \neq \emptyset$. 

   The proof in the presence of a horizon is exactly the same in the Cases (II) and (III) above, but 
the proof of Case (I) does not apply directly in this situation since $\cH$ gives rise to new boundary 
components of $M$. 

  Suppose first the distance of $\cH$ to $\dm$ is bounded below:
\be \label{a}
dist_{g_{i}}(\cH, \dm) \geq t_{0} > 0,
\ee
for all $i$. The apriori estimate \eqref{3.6} and the proof in Case (I) above apply and show 
that 
\be \label{aa}
|Rm| \leq k, \ \ |d\log u| \leq k,
\ee
hold in a neighborhood of fixed size $\sim t_{0}$ about $\dm$. However, $|Rm|$ may still blow up 
on approach to $\cH$. To prove that this is in fact not the case, we work on the Ricci-flat 4-manifold 
$N$, which has the single boundary component $\dn = \dm \times S^{1}$. The Chern-Gauss-Bonnet 
theorem for manifolds with boundary in this case gives 
$$\int_{N}|Rm|^{2} = 8\pi^{2}\chi(N) + \int_{\dn}(A*Rm + A^{3}),$$
where $A$ and $Rm$ are the $2^{\rm nd}$ fundamental form and curvature of $(N, g_{N})$ and $A*Rm$ 
and $A^{3}$ are algebraic expressions in $A$ and $Rm$. These are controlled by the corresponding 
3-dimensional data on $(M, g)$, together with $u$, $N(u)$. 

  By \eqref{aa}, the boundary term is bounded, and hence 
$$\int_{N}|Rm|^{2} = 8\pi^{2}\chi(N) + C.$$
Thus, $(N, g)$ is a Ricci-flat 4-manifold with a uniform bound on the $L^{2}$ of curvature. The same 
argument as in Case (A) gives a lower volume bound on $(N, g_{N})$. One also needs an upper diameter 
bound on $(N, g_{N})$ or equivalently $(M, g)$. The proof of Lemma 3.3 however does not apply 
directly in this situation, since $u$ is not bounded away from zero. Nevertheless, we claim 
that \eqref{3.14} still holds in the black hole case. Namely as discussed in the beginning of Section 2, 
$N(u) = \k = const > 0$ along any component of the horizon $\cH$, so that, by the divergence 
theorem applied to $\D u$, 
$$\int_{\dm}N(u) = c_{0} > 0.$$
Working then in the context of the proof of Lemma 3.3, one has then again by the divergence 
theorem
$$\int_{S(r)}N(u) = c_{0} > 0,$$
for all $r$. On $\hat M$, the estimate \eqref{3.6} gives $|N(u)| \leq \frac{k}{r}$ on $S(r)$ (since $u$ 
is bounded above) and hence
$$area(S(r)) \geq c_{1}r,$$
for some constant $c_{1} > 0$. This gives \eqref{3.14}. The rest of the proof of Lemma 3.3 proceeds 
as before, and establishes the upper diameter bound in the black hole case. 

   It follows from the compactness theorem in \cite{An0}, \cite{BKN}, that a subsequence of 
$(N, (g_{N})_{i})$ converges in the Gromov-Hausdorff topology to a smooth Ricci-flat orbifold 
$X$ with a finite number orbifold singularities of the form $C(S^{3}/\Gamma)$, (a cone on a 
spherical space form). However, such singularities are never static, i.e.~the metric near the 
orbifold singularity is not of the form \eqref{1.3}. Hence there are no singularities so that 
$X = N$ and the convergence of the metrics is smooth everywhere. 
 
   Next, to prove \eqref{a} does in fact hold, suppose not. Then (in a subsequence) $dist_{g_{i}}(\cH, \dm) 
= d_{i} \to 0$. Analogous to \eqref{3.7}, we rescale the metrics $(M, g_{i})$ by $d_{i}^{-2}$ so that 
$dist_{g_{i}'}(\cH, \dm) = 1$.  One may apply the same arguments as above to obtain smooth 
convergence everywhere to a limit. The limit has outer boundary flat $\bR^{2}$ with $H = 0$, so by the 
Case (I) argument again, the limit is flat, with $u = const$, say $u = 1$. However at distance 1, one has 
the blow-up limit of the horizon $\cH$ where $u = 0$. Since the convergence is smooth everywhere, 
this gives a contradiction.   

   The remaining parts of the proof in Case (I) carry over without further change. 
 
{\endproof}

  A standard consequence of Theorem 1.1 is the following result on the structure of the boundary 
map $\Pi$ on the non-degenerate metrics $\cE^{nd}$, (with respect to $\Pi$). 

\begin{corollary} \label{c3.4}
On any connected component $\cC_{+}^{nd}$ of $\cE_{+}^{nd}$, the boundary map 
$$\Pi: \cC_{+}^{nd} \to Met_{+}(\dm)\times \cC_{+} \times \bR^{+}$$
is a finite sheeted covering map onto its image. 
\end{corollary}

{\bf Proof:} The results above show that $\Pi$ is a local diffeomorphism and proper on $\cC^{nd}$. 
It is then standard that $\Pi$ is a covering map, necessarily finite sheeted (see for instance \cite{AP}). 

{\endproof}

\begin{remark}
{\rm (I) With the single exception of Lemma 3.2, most all of the results above apply with only minor 
changes to the exterior static vacuum problem - the static extension problem, cf.~\cite{AK}. Lemma 
3.2 is the main reason we discuss the interior problem in this paper; we hope to discuss the 
exterior problem further elsewhere. 

  (II) The assumption $K_{\g} > 0$ is used only in Cases (II) and (III) above, and it is not used very strongly 
there. We expect that Theorem 1.1 in fact holds for general boundary metrics $\g$, so that the condition 
$K_{\g} > 0$ is not necessary. It would be interesting to expore this further. 
}
\end{remark}

\begin{remark}
{\rm One may also consider static vacuum metrics with cosmological constant $\L$. These are 
given by solutions $(M, g, u)$ to the system
$$uRic_{g} = D^{2}u + \l u g, \ \ \D u = -\l u.$$
Again most all of the results above apply in this situation as well, with the exception of Lemma 3.3 
which will not hold in this generality when $\l < 0$. 

}
\end{remark}

\section{Degree Computation}

   As discussed in the Introduction, a smooth, proper Fredholm map $F: X \to Y$ of index 0 between 
Banach manifolds has a well-defined $\bZ_{2}$-valued degree, given by \eqref{1.9}. By Theorem 1.1, 
this applies to the map $\Pi$. We will also use the following simple generalization of \eqref{1.9}. 
Recall from Section 2 that a non-degenerate critical manifold of $\Pi$ is a connected $k$-dimensional 
submanifold $\S$ of $\cE$ such that the nullity of $(g, u)$ equals $k$, for each $(g, u) \in \cE$. 

\begin{proposition} \label{p4.1}
   Suppose $\Pi^{-1}(\g, [H], \mu)$ is a union of compact non-degenerate critical manifolds 
$\S_{i}$, $1 \leq i \leq m$. Then
\begin{equation} \label{4.1}
deg_{{\bZ_{2}}}(\Pi) = \sum_{i=1}^{m}\chi(\S_{i}) \ \ \   (mod \ 2). 
\end{equation}
\end{proposition}

{\bf Proof:} We may assume $m = 1$ with $\Pi^{-1}(\g, [H], \mu) = \S$. Choose a regular value 
$(\g', [H'], \mu')$ near $(\g, [H], \mu)$. By Theorem 1.1, $\Pi^{-1}(\g', [H'], \mu')$ is a 
finite set of points $P = \{p_{1}, \cdots, p_{r}\}$, $p_{i} = (g_{i}, u_{i})$. By Proposition 2.4, 
these points correspond bijectively to the critical points $q_{i}$ of a Morse function 
$\f: \S \to \bR$. (This is discussed in more detail in Proposition 5.2 below). By elementary Morse theory, 
the Euler characteristic $\chi(\S)$ (mod 2) equals the number of critical points of $\f$, 
(mod 2). 

{\endproof}

\begin{remark}
{\rm For mappings associated to variational problems, the $\bZ_{2}$-valued degree can sometimes 
be enhanced to a $\bZ$-valued degree. The origins of this are in the Morse theory for closed geodesics, 
as well as the Leray-Schauder degree; cf.~\cite{Bo}, \cite{Ni} for instance for surveys. Such a degree theory 
has also been developed for minimal surfaces by several authors, cf.~\cite{BT}, \cite{Tr}, \cite{W1}, \cite{W2}. 

  The main issue in defining a $\bZ$-valued degree is to construct  an orientation on the domain 
$X$ and target $Y$. For suitable variational problems the degree can then be calculated (as in the 
finite dimensional case) as 
\be \label{4.2}
deg F = \sum_{F(x_{i}) = y}(-1)^{ind(x_{i})},
\ee
where $ind(x_{i})$ is the index of the solution $x_{i}$, i.e.~the maximal dimension on which the 
$2^{\rm nd}$ variation is negative definite. 

   In the context of the proper Fredholm map $\Pi$ with the associated variational problem $\cL$ or 
$\w \cL$, the $2^{\rm nd}$ variation is given by the linearized operator $L$ in \eqref{2.21}, 
acting on the space $\cT_{0}^{m,\a}$ of deformations $h^{N}$  with fixed boundary conditions 
\eqref{2.20} at $\partial M$. Since $L$ acting on $\cT_{0}^{m,\a}$ is elliptic, Fredholm and self-adjoint, 
$L$ has a complete basis of eigenvectors in $L^{2}$, with discrete eigenvalues. From the form of 
\eqref{2.21}, one sees that $L$ is essentially a positive operator, i.e.~is bounded below, on 
transverse-tracefree deformations, so there are eigenvalues $\l_{i} \to + \infty$. However, on 
``pure trace forms", i.e.~divergence-free forms of the type $h^{N} = fg_{N} + \d^{*}Z$ with 
$Z = 0$ on $\dm$, $L$ is a negative operator. Thus, the index of $\cL$ (or $\w \cL$) is 
always infinite. This behavior is well-known in the context of closed manifolds, cf.~\cite{Be} for instance. 
It may be possible to define and compute a degree as in \eqref{4.2} using the spectral flow of the operator 
$L$ in place of the index. However, we will not pursue this further here. 

}
\end{remark}

  We now turn to the computation of the degree of $\Pi$ in some important cases. 

\medskip 

    First consider the case $\cH = \emptyset$, with $M = B^{3}$. Recall that $\mu$ is defined 
as in \eqref{1.7} as
$$\mu = \int_{\dm}(|d\nu|^{4} + \nu^{2})dv_{\g}.$$
Let $\g_{+1}$ denote the round metric of radius $1$ on $S^{2} = S^{2}(1)$ and $[2] = [c]$ the 
equivalence class of constant $H$. Let 
$$\S_{\mu} = \Pi^{-1}(\g_{+1}, [2], \mu) \subset \cE.$$ 
    
\begin{proposition}\label{p4.3}
For $M = B^{3}$, the inverse image $\S_{\mu}$ consists of the flat metric $g_{Eucl}$ on 
$B^{3}(1) \subset \bR^{3}$ with potential $u$ satisfying
\be \label{4.3}
u = 1 \ \ {\rm when} \ \ \mu = 0,
\ee
\be \label{4.4}
u = a + bz \ \ {\rm when} \ \ \mu > 0,
\ee
with $a + bz$ a general affine function on $B^{3}(1) \subset \bR^{3}$ with $|z| = 1$. Thus the 
inverse image $\S_{0}$ is a (critical) point, while the inverse image 
$\S_{\mu}$, $\mu > 0$ is a manifold diffeomorphic to $S^{2}\times S^{1}$.
\end{proposition}    

{\bf Proof:} Consider the set of static vacuum solutions $(g, u)$ with boundary data $\g = \g_{+1}$, 
$H = const$ and with fixed $\mu$. We claim that the metric $g$ is the flat metric $g_{Eucl}$ on 
the ball $B^{3}(1)$ of radius 1 and $u$ is an affine function $D^{2}u = 0$. To prove this, first 
by the result of Shi-Tam \eqref{3.10}, one has 
$$H \leq 2$$
with equality if and only if $(M, g)$ is flat. By the scalar constraint \eqref{2.1},
\be \label{4.4a}
u(|A|^{2} - H^{2} + 2) = 2(\D u + N(u)H),
\ee
so that
\be \label{4.5}
\int u[|A|^{2} - H^{2} + 2] = 0,
\ee
since $H = const$ and $\int_{\dm}N(u) = \int_{M}\D u = 0$. Next, let $A_{0}$ be the trace-free 
part of $A$. Then since $H^{2} \leq 4$, $|A|^{2} - H^{2} + 2 = |A_{0}|^{2} - 
\frac{1}{2}H^{2} + 2 \geq 0$. It follows from \eqref{4.5} that equality holds, and hence $|A|^{2} = 2$, 
so that $A = \g$ and $H = 2$. Thus, the Cauchy data $(\g, A)$ of $\dm$ in $M$ equal that of 
the round sphere $S^{2}(1) \subset \bR^{3}$. Regarding the Cauchy data for $u$ at $\dm$, 
the divergence constraint \eqref{2.2} implies that $N(u) = u - c = \bar u$ on $\dm$, where $c$ is the 
mean value of $u$ on $\dm$. From this and the scalar constraint \eqref{2.1} one obtains 
$\D \bar u + 2\bar u = 0$, so that $u$ is an eigenfunction of the Laplacian on $S^{2}(1)$. 
Hence, the Cauchy data for $u$ are that of an affine function on $\bR^{3}$, restricted to 
$S^{2}(1)$. Since all the Cauchy data are standard flat data, it follows from the unique 
continuation theorem in \cite{AH} (or from analyticity) that the solution itself is a standard flat 
solution $g_{Eucl}$, $u = a + bz$. 

   When $\mu = 0$, so $u = 1$, the variation $(k, u') = (0, 1)$ is in $Ker D\Pi$, so that 
$(M, g_{Eucl}, 1)$ is a critical point of $\Pi$. Next consider the family $(g_{Eucl}, u)$ with $u = 
a + bz$ with $\mu = \mu_{0} > 0$. The function $z$ is uniquely determined by a unit vector 
$z \in S^{2}(1)$. Moreover, the space of such potentials with $\mu$ fixed is compact and hence 
$a, b$ vary over a circle $S^{1}$ (topologically); the size of the circle is determined by $\mu$. 
(This can also be verified by direct computation). As $\mu\to 0$, $b \to 0$ and $a \to 1$. 
Hence $\S_{\mu} = \Pi^{-1}(\g_{+1}, [2], \mu) \simeq S^{2}\times S^{1}$. 

{\endproof}

   Observe that
\be \label{4.5a}
\S \equiv \cup_{\mu \geq 0}\S_{\mu} \simeq S^{2}\times \bR^{2},
\ee
where the origin in $\bR^{2}$ corresponds to $u = 1$, (so $a = 1$, $b = 0$). Note also 
that $\S = (\Pi_{B})^{-1}(\g_{+1}, 2)$. 

\begin{remark}\label{r4.4}
{\rm It is not unreasonable to conjecture that Proposition 4.3 holds for any flat boundary data, 
i.e.~any $(\g, H)$ which arise from an isometric embedding or immersion $\iota: (M, g_{flat}) 
\looparrowright \bR^{3}$. We note that $(\g, H)$ uniquely determine the isometric immersion $\iota$ 
into $\bR^{3}$ up to rigid motion. This follows from the well-known result (cf.~\cite{Sa} for instance) 
that there are no non-trivial Bonnet pairs for spheres $S^{2}$ immersed in $\bR^{3}$, i.e.~any two 
isometric immersions $S^{2} \looparrowright \bR^{3}$ with the same mean curvature differ by a 
rigid motion. (The difference of the two second fundamental forms is a Hopf differential on $S^{2}$ 
and hence vanishes). Call such a pair $(\g, H)$ flat boundary data (embedded flat boundary data 
if $\iota$ above is an embedding) and set 
\be \label{4.6}
\S_{\mu}(\g, H) = \Pi^{-1}(\g, [H], \mu).
\ee

  Thus one would like to know if for flat boundary data $(\g, H)$ solutions in $\Pi^{-1}(\g, [H], \mu)$ 
are uniquely realized by a flat metric $g_{flat}$ on $M$ with $u$ an affine function. The main issue in 
proving this is knowing the precise value of the mean curvature $H_{M} \in [H]$ of $(M, g)$. One has 
$H_{M} = \l H$ and by \eqref{3.10}, $\l \leq 1$. If $K_{\g} > 0$ and $\l = 1$ then the rigidity 
statement in \eqref{3.10} implies that $(M, g)$ is flat. Thus the issue is knowing if there are any 
static vacuum solutions with flat boundary data and with $\l < 1$, (and also dealing with the case 
of non-convex embeddings or immersions of $S^{2}$ into $\bR^{3}$). 

}
\end{remark}

   Next we claim $\S_{\mu}(\g, H)$ for flat boundary data $(\g, H)$ close to $(\g_{+1}, 2)$ and $\mu > 0$ 
is a non-degenerate critical manifold, i.e.~the nullity of any $(g_{flat}, u) \in \S_{\mu}$ equals 
$dim(S^{2}\times S^{1}) = 3$.

\begin{proposition}\label{p4.5}
For any flat boundary data $(\g, H)$ (sufficiently) near $(\g_{+1}, 2)$, the manifolds $\S_{\mu}(\g, H)$, 
$\mu > 0$, are non-degenerate critical manifolds for $\Pi$ in $\cE$. 
\end{proposition}

{\bf Proof:} One needs to show that if $(k, u')$ is an infinitesimal static Einstein deformation of 
$(M, g_{flat}, u)$ with $(k^{T}, [H'_{k}], \mu') = (0,0,0)$ then $(k, u') = (0, u')$ with $u'$ an affine 
variation of the affine function $u$. 

   Consider first the case $(\g, H) = (\g_{+1}, 2)$, so $(M, g)$ is the standard round unit ball 
$B^{3}(1) \subset \bR^{3}$. We assume for the moment $u = 1$, (so $\mu = 0$) and prove that 
the nullity of $(B^{3}(1), g_{flat}, 1)$ equals 4. 

   Since $u = 1$, the linearization of the scalar constraint \eqref{4.4a} in the direction $k$ at 
$\g_{+1}$, $H = 2$ gives
\be \label{4.7}
2(\D u' + 2N(u')) = 2\<A'_{k}, A\> - 2HH'_{k} = 2H'_{k} - 4H'_{k} = -4\l,
\ee
on $\dm$ where $H'_{k} = \l H = 2\l$. Here we have used the fact that $A = \g$ and $tr A'_{k} = 
H'_{k}$ since $k^{T} = 0$. Integrating \eqref{4.7} over $\dm$ gives 
$$-2\l area(\dm) = \int_{\dm}N(u') = \int_{M}\D u'.$$
Since $0 = (\D u)' = \D' u + \D u' = \D u'$ (again using $u = 1$), it follows that $\l = 0$ and hence 
\be \label{4.8}
\D u' + 2N(u') = 0.
\ee

  Now $\D u' = 0$ on $M = B^{3}(1)$ and at the boundary \eqref{4.8} holds, giving a relation between 
the Dirichlet and Neumann boundary data. Decompose $u'$ in terms of eigenfunctions of the round 
Laplacian on $S^{2}(1)$, (these are the restrictions of harmonic polynomials of degree $k$ on 
$\bR^{3}$ to $S^{2}(1)$), so that 
$$u' = \sum c_{k}\f_{k},$$
where $\D \f_{k} = -\l_{k}\f_{k}$, with $\l_{k} = k(k+1)$ (the eigenvalues of $\g_{+1}$). Hence  
$$\D u' = -\sum c_{k}k(k+1)\f_{k}.$$
Next, one has the Dirichlet-to-Neumann map $u' \to N(u')$. It is well-known, cf.~\cite{Ta} for instance, 
that this is a non-negative, elliptic and self-adjoint, first order pseudo-differential operator, with spectrum 
$\mu_{m} = m$. Moreover, as above, the eigenfunctions are restrictions of harmonic polynomials of 
order $m$ to $S^{2}(1)$. Hence
$$N(u') = \sum m c_{m}\f_{m}.$$
and \eqref{4.8} gives
$$ -\sum c_{k}k(k+1)\f_{k} + 2\sum m c_{m}\f_{m} = 0.$$
Pairing this with any $\f_{\ell}$ and integrating, it follows that
$$-\ell(\ell+1) + 2\ell = 0,$$
so $\ell = 1$ is the only solution. Thus the only solutions of \eqref{4.8} as boundary data of harmonic 
functions on $S^{2}(1)$ are the first eigenfunctions of the Laplacian, i.e.~restrictions of affine functions 
to $S^{2}(1)$. It follows that $u'$ is an affine function, $u' = a' + b'z'$, on $M$. 

  By the linearized static vacuum equations, one has then $u'Ric + uRic' = (D^{2})'u + D^{2}u'$ which 
thus gives $Ric_{k}' = 0$, so that $k$ is an infinitesimal flat deformation of $(B^{3}(1), g_{Eucl}, u)$ 
with $k^{T} = 0$ on $\dm$. It is well-known that convex surfaces in $\bR^{3}$ are 
infinitesimally rigid, and hence it follows that $k = 0$ (in divergence-free gauge). 

  Writing $u' = a' + b'z'$, one sees that the dimension of the space of affine variations $u'$ is $4$, so that 
the nullity of $(g_{flat}, 1)$ is $4$. 

   Now of course in the situation above, $\mu = 0$. Consider the nullity of $(M, g_{flat}, u)$ for $u$ close to $1$ 
with $\mu > 0$. By general principles, the nullity of any solution $(M, g, u)$ near $(M, g_{flat}, 1)$ is at most $4$. 
However, the $2^{\rm nd}$ variation of $\mu$ at $(g_{flat}, 1)$ is strictly positive: 
$$D^{2}\mu(u', u') > 0,$$
for any $u'$. Hence the nullity of any solution $(M, g, u)$ near $(M, g_{flat}, 1)$ is in fact at most $3$. 
The nullity thus equals 3 for all flat solutions $(M, g_{flat}, u)$, with $u = a + bz$, 
with $a \sim 1$, $b \neq 0$, $b \sim 0$, any $z$ with $|z| = 1$. This proves Proposition 4.5 for $\mu$ 
close to $0$. 

   With further work, the computations above can be extended to cover the case of general $\mu$. Instead, we 
give a more conceptual or geometric proof. 

   Take any $(M, g_{flat})$ with general $u = a + bz$. In the universal cover $\w N$, the flat 4-manifold has 
the metric form
$$g_{\w N} = (a + bz)^{2}dt^{2} + g_{flat}.$$
The domain $\w N$ is embedded or immersed in $\bR^{4}$. We pull back $\w g_{N}$ by the 
diffeomorphism $\psi_{c}: \bR^{4} \to \bR^{4}$, $\psi_{c}(x,y,z,t) = (x, y, z-c,t)$. Then 
the pullback 4-metric has the form
$$(a + b(z-c))^{2}dt^{2} + g_{flat}.$$
The coefficient of $dt^{2}$ (corresponding to $u^{2}$) may be multiplied by an arbitrary constant $k^{2}$. 
(This may in fact also be accomplished in the universal cover by the dilation $t \to kt$). This gives the metric 
$$\hat g_{\w N} = k^{2}(a + b(z-c))^{2}dt^{2} + g_{flat}.$$
One has $\hat u = k(a + b(z-c)) = k(a - bc) + kbz$, so $\hat b = kb$ which may now be arbitrary, choosing 
$k$ arbitrary (large). Similarly, $\hat a = k(1 - bc)$, so that choosing $c = b^{-1}(1 - k^{-1})$ gives 
$\hat a = 1$.  

  This shows that one can pull back the solution $(g_{flat}, a + bz)$ with $b$ arbitrary by a diffeomorphism, 
to obtain a solution $(\hat g_{flat}, \hat a + \hat b z)$ with $b$ small. If $(g_{t}, u_{t})$ is a static Einstein 
deformation of $(g_{flat}, a + bz)$ or an infinitesimal static Einstein deformation, then the pullback is also 
a static Einstein deformation of $(\hat g_{flat}, \hat a + \hat b z)$. Also, the boundary conditions are 
preserved by such a pullback. Since the nullity of $(\hat g_{flat}, \hat a + \hat b z)$ equals $3$, the same 
holds for the nullity of $(g_{flat}, a + bz)$. 

{\endproof}

   As in Remark 4.3, it would be very interesting to know if Proposition 4.5 holds for general flat domains 
$M \looparrowright \bR^{3}$. We conjecture this is in fact the case; the proof above however does not 
seem to generalize easily to this situation. 

\medskip 
\noindent
{\bf Proof of Theorem 1.2.}

\medskip 

   Since $\chi(\S_{\mu}) = 0$, the result follows from Propositions 4.1, 4.3 and 4.5. 
   
{\endproof}

\begin{remark}
{\rm There are other ways to see that $deg_{\bZ_{2}}\Pi = 0$ besides the proof above (although the 
proof above provides much further information that will be used in Section 5). 

    For example, the data $(\g, [H], 0)$ is not in $Im(\Pi)$, for any non-flat boundary data $(\g, H)$, 
i.e.~any $(\g, H)$ not arising from an isometric immersion $M \looparrowright \bR^{3}$. Namely  
if $\mu = 0$ then $u \equiv 1$ and so any static vacuum solution is flat. Since $\Pi$ is 
thus not surjective, one must have $deg_{\bZ_{2}}\Pi = 0$. 

   In a related vein, consider the behavior of $\Pi$ under rescalings $u \to c u$ of the potential $u$. 
This leaves the data $(\g, [H])$ invariant while $\mu$ becomes 
\be \label{mu}
\mu(d) = \mu + 2d\int_{\dm}\nu dv_{\g} + d^{2}\int_{\dm}1dv_{\g},
\ee
where $d = \log c$. The function $\mu(d)$ has a single critical point, a minimum at $d_{0} = 
\oint_{\dm}\nu dv_{g}$, where $\oint$ is the average value, with $\mu(d) \to +\infty$ as $d \to 
\pm \infty$. The map $\mu$ is thus a 2-1 map with a simple fold singularity at $d_{0}$. An easy 
computation using the H\"older inequality shows that $\mu(d_{0}) > 0$ unless $u \equiv 1$. 
Hence the inverse image $\Pi^{-1}(\g, [H], \mu)$ of any regular value always consists of an 
even number of points, with pairs $(g, u_{i})$, $i = 1, 2$, differing just by rescalings of the potential, 
whenever $\mu \neq \mu(d_{0})$. 

  One could naturally choose a normalization for the potential $u$, by choosing $c$ above so that 
$\w u = cu$ realizes the (unique) minimum of $\mu(d)$ above. Then $\cE_{norm} = \{(g, u)\} \in \cE$ 
such that $u = \w u$ is a codimension one smooth Banach manifold, and the induced map $\Pi_{[B]}$ 
from \eqref{1.6a}, 
$$\Pi_{[B]}: \cE_{norm} \to \cD,$$
$$\Pi_{[B]}(g, \w u) = (\g, [H]),$$ 
is smooth, Fredholm, of Fredholm index 0. However, $\Pi_{[B]}$ here is still not proper; (consider again the 
space of flat solutions with fixed $(\g, H)$ as discussed in the Introduction). 

   Observe that the arguments above hold for strict static solutions $M$ with arbitrary topology. 
}
\end{remark}

  Next we turn to the black hole case $\cH \neq \emptyset$. The simplest topology in this 
situation is $N = D^{2}\times S^{2}$ with $M = I\times S^{2}$ where $I = [0,1]$ (for instance) 
with the inner boundary $\{0\}\times S^{2}$ corresponding to the horizon $\cH$. The 
``standard solution" with this topology is the curve of Schwarzschild metrics
$$g_{Sch} = c^{2}V^{2}d\t^{2} + V^{-2}dr^{2} + r^{2}g_{S^{2}(1)},$$
where $V^{2} = (1 - \frac{2m}{r})$ and $\t \in [0, 2\pi]$, with $u = cV$. The smoothness of the 
metric at the horizon $\cH = \{r = 2m\}$ requires  
$$c = 4m,$$
so that
$$g_{Sch} = 16m^{2}(1-\frac{2m}{r})d\t^{2} + (1-\frac{2m}{r})^{-1}dr^{2} + r^{2}g_{S^{2}(1)}.$$
One has $r \geq 2m$ and we let $r \in [2m, R]$ with the outer boundary $\dm$ at the locus 
$\{r = R\}$. The boundary data $(\g, [H], \mu)$ are then of the form $(\g_{R}, [c], \mu)$ with 
$\g_{R} = R^{2}\g_{+1}$, $\g_{+1}$ the round metric of radius $1$ on $S^{2}$.  

   For convenience, set $R = 1$; (this may always be achieved by rescaling). Then on $\dm$, 
$\g = \g_{+1}$ and 
$$u^{2} = 16m^{2}(1 - 2m).$$
A simple calculation shows that (at $r = 1$),
$$H = 2\sqrt{1 - 2m}.$$ 
One has $m \in (0, \frac{1}{2})$ with $u \to 0$ on $\dm$ as $m \to 0$ or $m \to \frac{1}{2}$. The function 
$u = u(m)$ has a single maximum at $\frac{1}{3}$, corresponding to the well-known photon 
sphere $\{r = 3m\}$ of the Schwarzschild metric. For $m \neq \frac{1}{3}$, there are exactly two values $m_{\pm}$ 
of $m$ giving the same value for $u$ on $\dm$ and hence the same value to $\mu$. Thus on the 
Schwarzschild curve, the map $\Pi$ is two-to-one; 
$$\Pi^{-1}(g_{+1}, [c], \mu) \cap \{g_{Sch}\} = g_{Sch}(m_{\pm}) = (g_{\pm}, u_{\pm}).$$
These two branches (the small black hole for $m_{-} < \frac{1}{3}$ and the large black hole for 
$m_{+} > \frac{1}{3}$) merge together at $m = \frac{1}{3}$. The solution with $m = \frac{1}{3}$ is 
a critical point of $\Pi$. 
 
  This behavior on the curve of Schwarzschild metrics has long been known in the physics community, 
cf.~\cite{AG} for further discussion and examples. 

\begin{remark}
{\rm On the Schwarzschild curve, one may vary the mean curvature constant $H$, while keeping 
the boundary metric $\g$ fixed. The map $\Pi_{B}$ is not transverse to the fibers $[H]$, compare 
with Propositions 2.5 and 2.6. 
}
\end{remark}

  The discussion above suggests that 
$$deg_{\bZ_{2}} \Pi = 0,$$
for $N = D^{2}\times S^{2}$. This would follow from the following two conjectures.

\medskip 
\noindent
{\bf Conjecture:} (Global Uniqueness) For $N = D^{2}\times S^{2}$ with $\cH = S^{2}$, one has
$$\Pi^{-1}(\g_{+1}, [c], \cdot) = g_{Sch}(m_{\pm}).$$

\medskip 
\noindent
{\bf Conjecture:} (Infinitesimal Uniqueness)  The metrics $g_{Sch}(m)$ with $R = 1$ are regular points 
of $\Pi$, provided $m \neq \frac{1}{3}$. 
\medskip

  These conjectures are versions of the well-known black hole uniqueness theorem for the 
Schwarz-schild metrics, but with finite boundary. We conjecture similar statements hold for 
the boundary map $\Pi_{B}$ with boundary data $(\g, H)$. 

  Both conjectures above seem rather difficult to prove (in contrast to the analogous case when 
$\cH = \emptyset$) and it would be interesting to make progress on them. We make simply the 
following remark. 

\begin{proposition}\label{p4.7}
Suppose $(M, g, u) \in \Pi^{-1}(\g_{+1}, [c], \mu_{0})$, for some $\mu_{0} > 0$. If $(M, g, u)$ is 
a regular point of $\Pi$, then $(M, g, u)$ is the Schwarzschild metric $g_{Sch}(m)$ for some $m$. 
\end{proposition}

{\bf Proof:} By assumption, the derivative $D\Pi$ of the map $\Pi: \cE \to \cD \times \bR$ is an 
isomorphism at $(M, g, u)$. Let $X$ be any Killing field of $\g_{+1}$ and extend $X$ to a smooth 
vector field on $M$. By construction, the form $\k = (\cL_{X}g, \cL_{X}u) = (2\d^{*}X, X(u))$ is 
in the kernel $K$ of $D\Pi$. It follows that there is a vector field $Z$ on $M$ with $Z = 0$ on 
$\dm$ (so $Z$ is tangent to ${\rm Diff}_{1}(M)$) such that $(\cL_{X}g, \cL_{X}u) = (\cL_{Z}g, \cL_{Z}u)$. 
Hence $\w X = X - Z$ is a Killing field of $(M, g)$ preserving the potential $u$. 

   Since the isometry group of $(S^{2}(1), [c])$ is $SO(3)$, this extends to a group of isometries of 
$(M, g, u)$ so that $(M, g, u)$ is spherically symmetric. The result then follows easily. 

{\endproof}

\section{Morse theory on $\cE$}

  In this section, we use Lyusternik-Schnirelman and Morse theory methods to study local properties and 
bifurcation phenomena in the space $\cE$ with respect to the boundary maps $\Pi$ and $\Pi_{B}$ and 
prove Theorem 1.3.  

   In a neighborhood of a regular point $(g, u) \in \cE$, $\Pi$ (or $\Pi_{B}$) is a local diffeomorphism, so that 
nearby solutions are uniquely determined by the boundary data $(\g, [H], \mu)$, and similarly for 
$\Pi_{B}$. 

   Consider first the boundary map $\Pi_{B}$. Suppose $\S \subset \cE$ is a non-degenerate 
critical manifold with respect to $\Pi_{B}$. The structure of $\cE$ near $\S$ is described by Proposition 2.4. 
Each point $(g_{0}, u_{0})$ in $\S$ has a neighborhood $\cU$ in $\cE$ given as the locus 
$z^{-1}(0) \subset \cB \times K$, or more precisely given by $I(z^{-1}(0))$, where 
$$z(\g, H, \k) = E(I(\g, H, \k)).$$
The space $z^{-1}(0)$ is a smooth submanifold of $\cB\times K$ of codimension $k$. Fixing the 
boundary data of $\S$ to be $(\g_{0}, H_{0})$, one has $\Pi_{B}^{-1}(\g_{0}, H_{0}) \cap \cU 
= \S \cap \cU$, with tangent space $K$. Thus a neighborhood of $0 \in K$ gives rise to a local chart 
$\psi$ for $\S$ near $(g_{0}, u_{0})$. 

  Now fix (arbitrary) boundary data $(\g, H)$ near $(\g_{0}, H_{0})$, so that now $I: K \to Met^{m,\a}(M)\times 
C^{m,\a}(M)$. By Proposition 2.4, $\Pi_{B}^{-1}(\g, H) \cap \cU$ is exactly the set of critical 
points $(g, u)$ of $\cL$ in $Met^{m,\a}(M)\times C^{m,\a}(M)$ near $(g_{0}, u_{0}))$.  Thus, via the chart $\psi$, 
the functional $\cL \circ I \circ \psi$ is a function on a neighborhood $V$ of $(g_{0}, u_{0})$ in $\S$, 
\be \label{5.1}
\hat \cL \equiv \cL\circ I \circ \psi: V \subset \S \to \bR.
\ee
By construction, the nullity of the static Einstein metric $(g, u)$ in $Met^{m,\a}(M)\times C^{m,\a}(M)$ 
(with respect to $\Pi_{B}$) equals the nullity of the critical point of $\hat \cL$ in \eqref{5.1}. The 
same discussion holds with $\w \cL$ in place of $\cL$ in the black hole case $\cH \neq \emptyset$ 
and we let $\hat \cL$ denote the functional in both cases. 

   The local description above is canonical (independent of any coordinates) and smooth. If $\S$ is 
compact, so there only finitely many neighborhoods to patch together, one obtains the following:
   
\begin{proposition}\label{p5.1}
Suppose $\S$ is a compact non-degenerate critical manifold for $\Pi_{B}$. For fixed boundary data $(\g, H)$ 
near $\Pi_{B}(\S)$, the functional $\hat \cL: \S \to \bR$ has critical points exactly equal to the static Einstein 
metrics $(g, u) \in \cE$ near $\S$. Moreover, the nullity of $(g, u) \in \cE$ with respect to $\Pi_{B}$ equals 
the nullity of the corresponding critical point $\hat \cL: \S \to \bR$. 
\end{proposition}

{\endproof}

  Note however that Proposition 5.1 does not directly apply in the strict static case, since non-degenerate critical 
manifolds of $\Pi_{B}$ are always non-compact (since they are invariant under rescalings of the potential $u$). 
Thus we consider the analogous situation for the boundary map $\Pi$. 

   Let then $\S$ be a compact non-degenerate critical manifold of $\Pi$. By Theorem 1.1, $\S$ is automatically 
compact if $\Pi(\S) = (\g, [H], \mu)$ satisfies $K_{\g} > 0$. 

\begin{proposition}\label{p5.2}
In the strictly static case, for any fixed boundary data $(\g, [H], \mu)$ near $\Pi(\S)$, the functional 
$\hat \cL: \S \to \bR$ as above has critical points exactly equal to the static Einstein metrics 
$(g, u) \in \cE$ near $\S$. Again, the nullity of $(g, u) \in \cE$ with respect to $\Pi$ equals the 
nullity of the corresponding critical point $\hat \cL: \S \to \bR$. 
\end{proposition}

{\bf Proof:} The proof is basically the same as that of Proposition 5.1. One only needs to define the map 
$z$ suitably and this amounts to redefining the equation \eqref{2.23} so that it is adapted to the 
boundary map $\Pi$ in place of the previous $\Pi_{B}$. As before the linearization $L$ maps 
$$L: \hat \cT_{0}^{m,\a} \to \cT^{m-2,\a},$$
where $\hat \cT_{0}$ denotes divergence-free forms $h^{(N)}$ satisfying the $\Pi$-boundary condition
$$(h^{T}, [H]'_{h}, \mu'_{(h, u')}) = (0,0,0).$$
To compare with this with $\cT_{0}$ in \eqref{2.20}, consider the boundary data $(0, \l H, 0)$ with 
(say) $\l = 1$ and extend this to a divergence-free form $h_{0}^{N}$ on $N$. Let $L(h_{0}^{N}) = f_{0}$. 
Then $h_{0}^{N} \in \hat \cT_{0}$, but $h_{0}^{N} \notin \cT_{0}$. Also recall that any form 
$h^{N} = (h, u') \in \cT_{0}$ may be ``trivially" shifted by adding a term $(0, \l u)$ in $K_{\Pi_{B}}$ 
so that $h^{N} \in \hat \cT_{0}$, cf.~also Remark 4.6. This gives an isomorphism
$$\cT_{0}\oplus \<h_{0}^{N}\> \to \hat \cT_{0}.$$

   Next, let $K_{\Pi}^{\perp}$ be the $L^{2}$ orthogonal complement of $K_{\Pi}$ 
in $\hat \cT^{m-2,\a}(N)$ and let $K_{\Pi,0}^{\perp}$ be the $L^{2}$ orthogonal complement of $K_{\Pi}$ 
in $\hat \cT_{0}^{m,\a}(N)$. We claim that as in \eqref{2.23}, $L$ induces an isomorphism
\be \label{5.2}
\hat L: K_{\Pi,0}^{\perp} \to K_{\Pi}^{\perp}.
\ee
To see this, recall from Corollary 2.7 that if $k \in K_{\Pi}$ then $k \in K_{\Pi_{B}}$ in the strict static case, 
so $K_{\Pi} \subset K_{\Pi_{B}}$ and hence, by \eqref{2.38}, there is an isomorphism (inclusion) 
$$\iota: K_{\Pi_{B},0}^{\perp} \to \cW \subset K_{\Pi,0}^{\perp},$$
where $\cW$ is a codimension 1 hypersurface. This gives an isomorphism $\iota^{-1}: \cW \to 
K_{\Pi_{B},0}^{\perp}$ and also an isomorphism $\cW \oplus \<h_{0}^{N}\> \simeq K_{\Pi,0}^{\perp}$, 
(since $\cW \subset K_{\Pi_{B},0}^{\perp}$ and $h_{0}^{N} \notin  K_{\Pi_{B},0}^{\perp}$). Combining 
these gives an isomorphism 
$$\iota_{1}: K_{\Pi, 0}^{\perp} \to K_{\Pi_{B},0}^{\perp} \oplus \<h_{0}^{N}\>.$$
Applying then $L$ as in \eqref{2.23} gives the isomorphism 
$$L: K_{\Pi_{B},0}^{\perp} \oplus \<h_{0}^{N}\> \to K_{\Pi_{B}}^{\perp} \oplus \<f_{0}\>.$$
Composing this with isomorphism induced again by $\iota$ 
$$\iota_{2}: K_{\Pi_{B}}^{\perp} \oplus \<f_{0}\> \to K_{\Pi}^{\perp} \subset \hat \cT^{m-2,\a},$$
gives the claim \eqref{5.2} with $\hat L = \iota_{2}\circ L \circ \iota_{1}$. 

{\endproof}

   Propositions 5.1 and 5.2 allow one to use the tools of Lyusternik-Schnirelman theory and Morse theory to 
study the structure of static vacuum solutions with boundary data near $\Pi(\S)$ (or $\Pi_{B}(\S)$). Recall 
that the Lyusternik-Schnirelman theory implies that the number of critical points of a smooth function on 
a compact manifold $S$ is at least $Cat(S)$, where $Cat$ is the smallest number of open contractible sets 
covering $S$. Similarly, Morse theory implies that the number of critical points of a Morse function is at least 
the sum of the Betti numbers of $S$. 

\medskip 

  The results above lead easily to the proof of Theorem 1.3. 

\medskip
\noindent 
{\bf Proof of Theorem 1.3.}

   For any flat boundary data $(\g_{0}, H_{0})$ near $(\g_{+1}, 2)$, the manifolds $\S_{\mu}$ defined as in 
\eqref{4.6} with $\mu > 0$ are compact non-degenerate critical manifolds for $\Pi$, by Proposition 4.5. 
By Proposition 5.2, for any $(\g, [H], \mu)$ close to $(\g_{0}, [H_{0}], \mu)$, static vacuum solutions in 
$\Pi^{-1}(\g, [H], \mu)$ near $\S_{\mu}$ are in one-to-one correspondence with the critical points of a 
smooth function $\hat \cL: \S_{\mu} \to \bR$. One has $\S_{\mu} \simeq S^{2}\times S^{1}$ and 
$Cat(S^{2}\times S^{1}) = 3$. Thus for each $(\g, [H], \mu)$ near flat boundary data, there are at least 
3 distinct static vacuum solutions $(M, g_{i}, u_{i})$, $1 \leq i \leq 3$, with 
\be \label{5.9}
\Pi(g_{i}, u_{i}) = (\g, [H], \mu).
\ee
By the Sard-Smale theorem, the regular values of $\Pi$ are generic; in fact by Propositions 2.5 and 2.6 
the critical values of $\Pi$ are of codimension at least one in $\cD \times \bR^{+}$. Hence for 
generic boundary data $(\g, [H], \mu)$, the functional $\hat \cL: \S_{\mu} \to \bR$ is a Morse function 
on $\S_{\mu} \simeq S^{2}\times S^{1}$ and so has at least 4 critical points $q_{i}$, with $index(q_{i}) 
= i-1$ for each $1 \leq i \leq 4$, i.e.~a minimum, maximum and two saddle points. This gives at 
least 4 distinct solutions to \eqref{5.9}. 

{\endproof}

  Theorem 1.2 implies there is an even number of solutions in $\Pi^{-1}(\g, [H], \mu)$ for any regular 
value $(\g, [H], \mu)$ of $\Pi$. By Remark 4.6, pairs of solutions $(g, u_{i})$, $i = 1,2$ in 
$\Pi^{-1}(\g, [H], \mu)$ just differ by rescaling the potential $u$, $u_{2} = cu_{1}$ for some $c > 0$. 
These rescaling pairs are essentially geometrically equivalent (the $4$-manifolds $(N, g_{N})$ are locally 
isometric) and so should be considered as just one solution. Of course the rescaling pairs 
$(g, u_{i})$ merge together at the rescaling critical point $\mu = \mu_{0}$ as following \eqref{mu}. 
Thus Theorem 1.3 gives at least two geometrically distinct solutions in $\Pi^{-1}(\g, [H], \mu)$, 
$\mu > 0$, for any $(\g, [H])$ close to $(\g_{+1}, [2])$.

   Recall from Corollary 3.4 that the boundary map 
$$\Pi: \cC_{+}^{nd} \to Met_{+}(\dm)\times \cC_{+} \times \bR^{+}$$
is a finite sheeted covering map onto its image, where $\cC_{+}^{nd}$ is any connected component 
of $\cE_{+}^{nd}$. Consider first this map in a neighborhood of the (rather large) space $\cE_{flat,+}$ of 
flat solutions in $\cE_{+}$, so $K_{\g} > 0$. By the Weyl embedding theorem, $H$ is then uniquely 
determined by $\g$ so that there is a diffeomorphism
$$\cE_{flat,+} \simeq \cP \times S^{2}\times D^{2},$$ 
where $\cP$ is the space of metrics of positive Gauss curvature on $S^{2}$ and $S^{2}\times D^{2}$ 
represents the space of affine potential functions as in \eqref{4.5a}. By \cite{RS}, the space $\cP$ is contractible, 
so that $\cE_{flat,+}$ is homotopy equivalent to $S^{2}$ and hence in particular is simply connected. 

  Localizing now to a neighborhood of standard round boundary data $(\g_{+1}, [2], \mu)$, consider the 
$\e$-ball $B_{\e}(\g_{+1}, [2], \mu)$ about $(\g_{+1}, [2], \mu)$ in $Met_{+}(\dm)\times \cC_{+} \times \bR^{+}$. 
and let $\cW = B_{\e}(\g_{+1}, [2], \mu)\setminus \Pi(\cE_{flat})$ be the complement of flat boundary data 
in $B_{\e}(\g_{+1}, [2], \mu)$. It follows then from Theorem 1.3 (and Propositions 4.3 and 4.5) that for 
$\e > 0$ sufficiently small, 
$$\Pi^{-1}(\cW) = \O = \cup \O_{i} \subset \cE^{nd},$$ 
where $\O_{i}$ are the components of $\O$, $1 \leq i \leq 2k$, $k \geq 2$, and the induced maps on 
each component 
\be \label{5.10}
\Pi: \O_{i} \to \cW
\ee
are diffeomorphisms for each $i$. 

   To understand this in more detail, chooose any flat boundary data $(\g_{0}, H_{0})$ near 
$(\g_{+1}, 2)$ and choose any regular value $(\g, [H], \mu)$ of $\Pi$ sufficiently near 
$(\g_{0}, [H_{0}], \mu)$. Then $\Pi^{-1}(\g, [H], \mu) = \{(g_{i}, u_{i})\}$, $1 \leq i \leq 2k$. 

   Now consider the structure of $\Pi^{-1}(\g, [H], t\mu) \subset \cE$ for $t \in [0,1]$. Since $\Pi$ is proper, 
this is a compact set in $\cE$, given as an even number of curves $\s_{i}(t) = (g_{i}(t), u_{i}(t))$ off 
$\cE^{d}$. For $t$ close to 1, each curve $\s_{i}(t) \in \O_{i}$. Pairs of such curves are rescaling curves. 
However, at $t = 0$, 
$$\Pi^{-1}(\g, [H], 0) = \emptyset,$$
since any solution with $\mu = 0$ is necessarily flat and $(\g, [H])$ are not flat boundary data, 
by construction. Hence for $t$ sufficiently close to $0$ there are also no solutions in $\Pi^{-1}(\g, [H], t\mu)$. 
In particular, $\e$ above depends on $\mu$ with $\e \to 0$ as $\mu \to 0$. 

  Thus the curves $\s_{i}(t)$ cannot project via $\Pi$ to $(\g, [H], t\mu)$ for $t$ small. They must thus end, 
and presumably join, at the critical locus $\cE^{d}$. The rescaling pairs meet when $\mu = \mu(d_{0})$, the 
minimal value of $\mu$. This corresponds to the passage from (say) 4 critical points to 3 critical points by the 
merging of a distinct pair (as in the passage from the Morse description to the Lyusternik-Schnirelman 
description above). The locus of these points $(g, u)$, $\mu = \mu(d_{0})$ form components of $\cE^{d}$ 
which are of codimension 1 in $\cE$; they project to a variety of critical values $(\g, [H], \mu(d_{0}))$ of 
codimension 1 in the target space. Of course $\cE_{flat} \subset \cE^{d}$ is of infinite codimension in $\cE$. 

\medskip 

   One sees that even in a neighborhood of standard flat solutions, the structure of the space $\cE$ is 
surprisingly complicated. It would be of interest to explore this further.

\medskip 

   Finally, it is interesting to relate the discussion above with recent work of Jauregui \cite{J}, cf.~also \cite{JMT}. 
Thus, choose any $\g \in Met_{+}(S^{2})$, so $K_{\g} > 0$ and pick any positive function $H > 0$. 
For $\l$ sufficiently small, it is shown in \cite{J} that $(\g, \l H)$ bounds a metric on $M = B^{3}$ with 
nonnegative scalar curvature. On the other hand, by the Shi-Tam result \eqref{3.10}, 
this cannot be the case if $\l$ is too large, i.e. 
$$\l > \frac{\int_{\dm}H_{0}dv_{\g}}{\int_{\dm}Hdv_{\g}}.$$
It is also shown in \cite{J} that there is a unique $\l_{0} = \l_{0}(\g, H) > 0$ such that $(\g, \l H)$ has a 
nonnegative scalar curvature filling for $\l < \l_{0}$ and has no such filling for $\l > \l_{0}$. The conjecture is 
made that $(\g, \l_{0}H)$ has a filling which is static vacuum. 

   By Theorem 1.3 and the discussion above, for any $(\g, [H])$ close to flat round boundary data 
$(\g_{+1}, [2])$ there are at least 2 geometrically distinct static vacuum fillings $(g_{i}, u_{i})$ (not differing 
by rescalings of the potential) of the boundary data $(\g, [H])$. Another interesting open question is 
to determine whether $H \in [H]$ is the same or different for all such fillings and whether the 
transition value $\l_{0}H$ is realized by one (or all) of these static vacuum fillings.

\bibliographystyle{plain}

\begin{thebibliography}{WWW}

\footnotesize

\bibitem [1]{AG} M. Akbar and G. Gibbons, {\it Ricci-flat metrics with $U(1)$ action and the Dirichlet 
boundary-value problem in Riemannian quantum gravity and isoperimetric inequalities}, 
Class. Quantum Gravity, {\bf 20}, (2003), 1787-1822.  


\bibitem[2]{AP} A. Ambrosetti and G. Prodi, {\it A Primer of Nonlinear Analysis}, Cambridge 
Studies in Adv. Math, {\bf 34}, Cambridge Univ. Press, Cambridge, UK, (1993). 


\bibitem[3]{An0} M. Anderson, {\it Ricci curvature bounds and Einstein metrics on compact 
manifolds}, Jour. Amer. Math. Soc., {\bf 2}, (1989), 455-490. 

\bibitem[4]{An1} M. Anderson, {\it On stationary solutions to the vacuum Einstein equations},
Annales Henri Poincar\'e, {\bf 1}, (2000), 977-994.

\bibitem[5]{An1.5} M. Anderson, {\it On the structure of solutions to the static vacuum 
Einstein equations}, Annales Henri Poincar\'e, {\bf 1}, (2000), 995-1042. 

\bibitem[6]{An2} M. Anderson, {\it On boundary value problems for Einstein metrics},
Geom. \& Topology, {\bf 12}, (2008), 2009-2045.

\bibitem[7]{An3} M. Anderson, {\it Boundary value problems for metrics on $3$-manifolds}, in: 
Metric and Differential Geometry in Honor of J. Cheeger, Eds. X. Dai and X. Rong, Birkh\"auser Verlag, 
Basel, (2012), 3-17. 

\bibitem[8]{An4} M. Anderson, {\it Conformal immersions with prescribed mean curvature 
in $\bR^{3}$}, (preprint, 2012), arXiv:1204.5225 [math.DG] 


\bibitem[9]{AH} M. Anderson and M. Herzlich, {\it Unique continuation results for Ricci curvature
and applications}, Jour. Geom. \& Physics, {\bf 58}, (2008), 179-207; {\it Erratum}, ibid., {\bf 60},
(2010), 1062-1067.

\bibitem[10]{AT} M. Anderson, A. Katsuda, Y. Kurylev, M. Lassas and M. Taylor, {\it Boundary
regularity for the Ricci equation, geometric convergence and Gel'fand's inverse boundary
problem}, Inventiones Math., {\bf 158}, (2004), 261-321.

\bibitem [11]{AK} M. Anderson and M. Khuri, {\it On the Bartnik extension problem for the static 
vacuum Einstein equations}, Classical \& Quantum Gravity, (to appear (2013)), 
arXiv:0909.4550 [math.DG]. 

\bibitem[12]{BKN} S. Bando, A. Kasue and H. Nakajima, {\it On a construction of coordinates at 
infinity on manifolds with fast curvature decay and maximal volume growth}, Inventiones Math., 
{\bf 97}, (1989), 313-349. 


\bibitem[13]{Ba1} R. Bartnik, {\it New definition of quasi-local mass}, Phys. Rev. Lett., {\bf 62},
(1989), 2346-2348.


\bibitem[14]{Ba2} R. Bartnik, {\it Mass and 3-metrics of non-negative scalar curvature},
Proc. Int. Cong. Math., vol II, Beijing (2002), 231-240, Higher Ed. Press, Beijing, 2002.


\bibitem[15]{Be} A. Besse, {\it Einstein Manifolds}, Springer Verlag, Berlin, (1987).


\bibitem[16]{BT} R. B\"ohme and A. Tromba, {\it The index theorem for classical minimal surfaces}, 
Annals of Math., {\bf 113}, (1981), 447-499. 

\bibitem[17]{Bo} R. Bott, {\it Lectures on Morse theory, old and new}, Bulltein Amer. Math. Soc,. 
{\bf 7}, (1982), 331-358. 

\bibitem[18]{Ga} G. Galloway, {\it On the topology of black holes}, Comm. Math. Phys., {\bf 151}, 
(1993), 53-66. 

\bibitem[19]{GT} D. Gilbarg and N. Trudinger, {\it Elliptic Partial Differential Equations
of Second Order}, $2^{\rm nd}$ Edition, Springer Verlag, New York, (1983).

\bibitem[20]{J} J. Jauregui, {\it Fill-ins of nonnegative scalar curvature, static metrics and 
quasi-local mass}, Pacific J. Math., {\bf 261}, (2013), 417-444. 

\bibitem[21]{JMT} J. Jauregui, P. Miao \& L.-F. Tam, {\it Extensions and fill-ins with 
nonnegative scalar curvature}, arXiv: 1304.0721 [math.DG] 

\bibitem[22]{K} K. Knox, {\it A compactness theorem for Riemannian manifolds with boundary 
and applications}, (preprint, 2012), arXiv:1211.6210 [math.DG] 

\bibitem [23]{Mo} C. B. Morrey, Jr., {\it Multiple Integrals in the Calculus of Variations}, Springer 
Verlag, New York, (1966). 

\bibitem[24]{Ni} L. Nirenberg, {\it Variational and topological methods in nonlinear problems}, 
Bulletin Amer. Math. Soc., {\bf 4}, (1981), 267-302. 

\bibitem[25]{Pe} P. Petersen, {\it Riemannian Geometry}, $2^{\rm nd}$ Edition, Springer Verlag,
New York, (2006).

\bibitem[26]{RS} J. Rosenberg and S. Stolz, {\it Metrics of positive scalar curvature and connections with 
surgery}, in: Surveys on Surgery Theory, Annals of Math. Studies, Princeton Univ. Press, {\bf 149}, 
(2001), 353-386. 

\bibitem[27]{Sa} I. Kh. Sabitov, {\it Isometric surfaces with a common mean curvature and the problem of 
Bonnet pairs}, Sbornik: Mathematics, {\bf 203}, (2012), 111-152. 

\bibitem[28]{ST} Y.-G. Shi and L.-F. Tam, {\it Positive mass theorem and the boundary behaviors 
of compact manifolds with nonnegative scalar curvature}, Jour. Diff. Geom., {\bf 62}, (2002), 
79-125.  


\bibitem[29]{Sm} S. Smale, {\it An infinite dimensional version of Sard's theorem},
Amer. Jour. Math., {\bf 87}, (1965), 861-866.

\bibitem[30]{Ta} M. Taylor, {\it Partial Differential Equations II}, Applied Math. Sciences, {\bf 116}, Springer 
Verlag, New York, 1996. 

\bibitem[31]{To} P. Tod, {\it Spatial metrics which are static in many ways}, Gen. Rel. \& Gravitation, 
{\bf 32}, (2000), 2079-2090. 

\bibitem[32]{Tr} A. Tromba, {\it Degree theory on oriented infinite dimensional varieties and
the Morse number of minimal surfaces spanning a curve in ${\mathbb R}^{n}$}, Trans. Amer.
Math. Soc., {\bf 200}, (1985), 385-415.

\bibitem[33]{W1} B. White, {\it The space of $m$-dimensional surfaces that are stationary for
a parametric elliptic functional}, Ind. Univ. Math. Jour., {\bf 36}, (1987), 567-602.

\bibitem[34]{W2} B. White, {\it The space of minimal submanifolds for varying Riemannian
metrics}, Ind. Univ. Math. Jour., {\bf 40}, (1991), 161-200.


\end{thebibliography}

\end{document}